\documentclass[preprint,sort&compress,12pt]{elsarticle}
\usepackage{bbm}
\usepackage{amstext}
\usepackage{amsmath}
\usepackage{amssymb}
\usepackage{mathrsfs}
\usepackage{stmaryrd}
\usepackage{bm}
\usepackage{pstricks}
\usepackage{pst-coil}
\usepackage{pst-3d}
\usepackage{amsfonts}
\usepackage{amsthm} 
\usepackage{xcolor}

\newcommand{\mf}{\mathbf}
\newcommand{\mm}{\mathrm}
\newcommand{\ml}{\mathcal}

\newcommand{\be}{\begin{equation}}
\newcommand{\bea}{\begin{equation}\begin{aligned}}
\newcommand{\beas}{\begin{equation*}\begin{aligned}}
\newcommand{\eeas}{\end{aligned}\end{equation*}}
\newcommand{\eea}{\end{aligned}\end{equation}}
\newcommand{\ee}{\end{equation}}

\renewcommand{\div}{{\rm div }}

\begin{document} 
\begin{frontmatter}
\title{
Global Solutions of Three-dimensional Inviscid MHD Fluids with \\ Velocity Damping in Horizontally Periodic Domains
}

\author[sJ]{Fei Jiang}
\ead{jiangfei0591@163.com}
\author[FJ]{Song Jiang}
 \ead{jiang@iapcm.ac.cn}
\author[sJ]{Youyi Zhao \corref{cor1}}
\cortext[cor1]{Corresponding author. }
\ead{zhaoyouyi957@163.com}
\address[sJ]{College of Mathematics and
Computer Science, Fuzhou University, Fuzhou, 350108, China.}
 \address[FJ]{Institute of Applied Physics and Computational Mathematics,
 Beijing, 100088, China.}

\begin{abstract}
The \emph{two-dimensional} (2D) existence result  of global(-in-time) solutions for the motion equations of incompressible, inviscid, non-resistive magnetohydrodynamic (MHD) fluids with velocity damping had been established in [Wu--Wu--Xu, SIAM J. Math. Anal. 47 (2013), 2630--2656].  This paper further studies the existence of global solutions for the \emph{three-dimensional} (a dimension of real world) initial-boundary value problem in a horizontally periodic domain with finite height. Motivated by  the multi-layers energy method introduced in [Guo--Tice, Arch. Ration. Mech. Anal. 207 (2013), 459--531], we develop a new type of two-layer energy structure to overcome the difficulty arising from three-dimensional nonlinear terms in the MHD equations, and thus prove the initial-boundary value problem admits a unique global solution. Moreover the solution has the exponential decay-in-time  around some rest state. Our two-layer energy structure enjoys two  features: (1) the lower-order energy (functional)  can not be controlled by the higher-order energy. (2) under the \emph{a priori} smallness assumption of lower-order energy, we first close the higher-order energy estimates, and then further close the lower-energy estimates in turn.
\end{abstract}
\end{frontmatter}

\newtheorem{thm}{Theorem}[section]
\newtheorem{lem}{Lemma}[section]
\newtheorem{pro}{Proposition}[section]
\newtheorem{concl}{Conclusion}[section]
\newtheorem{cor}{Corollary}[section]
\newproof{pf}{Proof}
\newdefinition{rem}{Remark}[section]
\newtheorem{definition}{Definition}[section]

\section{Introduction}\label{introud}
\numberwithin{equation}{section}

In this paper, we investigate
the global(-in-time) well-posdedness of  the following \emph{three-dimensional} initial-boundary value problem for the incompressible, inviscid, non-resistive magnetohydrodynamic (MHD) fluid  with velocity damping in a horizontally periodic domain $\Omega$ with finite height:
\begin{equation}\label{0101}
\begin{cases}
\rho v_{t}+ \rho v\cdot\nabla v
+\nabla \left(P+\lambda|M|^2/2\right)+a \rho v =\lambda M\cdot\nabla M, \\
 {M}_{t}+ v\cdot\nabla {M}=M\cdot\nabla v, \\
\div  v =\mm{div}M=0,\\
(v,M)|_{t=0}=(v^0,M^0),\\
v_3|_{\partial\Omega } =0.
\end{cases}
\end{equation}
Next we shall explain the mathematical notations appearing in the problem above.

The unknowns ${v}:={v}(x,t)$, ${M}:={M}(x,t)$ and ${P}:={P}(x,t)$ denotes the velocity and magnetic field and  kinetic pressure of incompressible MHD fluids, resp., where $(x,t)\in \mathbb{R}^3\times \mathbb{R}_0^+$ and $\mathbb{R}_0^+:=[0,\infty)$.
The positive constants $\rho$, $a$ and $\lambda$ stands for the density of fluids,
  damping coefficient and  permeability of the vacuum, resp.
The horizontally periodic domain $\Omega$ with finite height is defined as follows:
\begin{align}\label{0101a}
\Omega:=\{(x_{\mm{h}}, x_3)^{\mm{T}}~|~x_{\mm{h}}:=(x_1, x_2)^{\mm{T}}
\in\mathbb{T},\;\;x_3\in(0,h)\}\subset\mathbb{R}^3,
\end{align}
where $\mathbb{T}=(l_1\mathbb{R}/\mathbb{Z})\times (l_2\mathbb{R}/\mathbb{Z})$,
 $\mathbb{R}/\mathbb{Z}$ represents the usual $1$-torus, and the superscript ${\mm{T}}$ denotes the transposition. Since $\Omega$ is horizontally periodic,
the set $\mathbb{T}\times \{0,h\}$, denoted by $\partial\Omega$, customarily represents the boundary of $\Omega$.
  Without loss of generality, we assume that $l_1=l_2=h=1$ for the sake of simplicity.

For the the viscous case, i.e. the equations \eqref{0101}$_1$--\eqref{0101}$_3$ with the viscous term  $-\nu \Delta v$  in place of the velocity damping term $a \rho v$,  the global  well-posedness, in which the initial data is a small perturbation around a non-zero trivial stationary state
{ (i.e., $(v,M)=(0,\bar{M})$, where $\bar{M}$ is a non-zero constant vector and $\bar{M}$ often called the impressive magnetic field)},
has been widely investigated, see \cite{lin2015global,zhang2014elementary,ZTGS} and \cite{ABIHZPOTG,xu2015global} for the 2D and 3D
 Cauchy problems resp., and see \cite{ren2016global} and \cite{TZWYJGw} for the 2D and 3D initial-boundary value problems resp.. The existence of global solutions to the 2D Cauchy problem with large initial perturbations was obtained by Zhang under a strong impressive magnetic field \cite{ZTGS}. {As for the well-poseness  of the 3D Cauchy and initial-boundary value problems  with large initial perturbations, to our best knowledge,  all available results are about the local(-in-time) existence, see \cite{fefferman2017local,HOCELETFCMDSFA,chemin2016local} for examples.
 We mention here that the corresponding compressible case has been also widely studied, see \cite{LYSYZG,LXLSNWDH,WJHWYF} and the references
 cited therein.}

 For the the inviscid case, Bardos--Sulem--Sulem used hyperbolicity of \eqref{0101}$_1$--\eqref{0101}$_3$ with $\mu=0$  to establish an interesting global existence result of classical
solutions of the Cauchy problem with small initial data in the H\"older space $H^s(\mathbb{R}^3)$ \cite{BCSCSPLL}; also see \cite{CYLZGW,HLBXLYP} for the case
of Sobolev spaces). We remark that such a result is not known
in the 3D incompressible Euler equations (i.e. the case of the absence of magnetic field).

For the inviscid case with a velocity damping term, Wu--Wu--Xu first given the existence of a unique global solution with algebraic decay-in-time for the 2D Cauchy problem of  \eqref{0101}$_1$--\eqref{0101}$_4$ with small initial perturbation  \cite{WJHWYFXXJG}. Recently, Du--Yang--Zhou also obtained the existence of a unique global solution with exponential decay-in-time for the initial-boundary value problem in a 2D slab domain with small initial perturbation around some non-trivial equilibrium \cite{DYYWZYOSJMA}.
To our best knowledge, there are not any available result for the 3D case,
 thought  the 3D well-posedness problem of inviscid fluids with velocity damping had been widely investigated, see \cite{TZWY2013JDE,PRHZK2009JDE,WWKYT2001JDE,STCTBWDH2012CPDE,HFMPRHARMA2003} for examples.  \emph{The goal of this paper is devoted to further providing the first 3D result by developing  a new type of two-layer energy structure to overcome the difficulty arising from  three-dimensional nonlinear terms in the MHD equations.}
It should be noted that the mathematical methods adopted in \cite{WJHWYFXXJG,DYYWZYOSJMA} for the well-posedness problem in
Eulerian coordinates cannot be applied to our 3D problem, therefore  next we shall reformulate the initial-boundary value problem \eqref{0101} in Lagrangian coordinates as in \cite{WYJ2019ARMA,JFJSSETEFP,JFJSOUI}.

\subsection{Reformulation in Lagrangian coordinates.}\label{subsec:02}

Let the flow map $\zeta$ to be the solution to the initial-value problem
\begin{equation}
\nonumber
\begin{cases}
  \zeta_t(y,t)=v(\zeta(y,t),t)&\mbox{in }\Omega,
\\
\zeta(y,0)=\zeta^0(y)&\mbox{in }\Omega,
\end{cases}
\end{equation}
where the
invertible mapping $\zeta^0:=\zeta^0(y):\Omega\to \Omega$ satisfies
\begin{align}
&\label{zeta0inta}
J^0:=\det(\nabla \zeta^0)=1  \mbox{ in } \Omega,\\
&\label{zeta0inta0}
\zeta^0_3 =y_3  \mbox{ on } \partial\Omega .
\end{align}
Here and in what follows  $\det$ denotes a determinant of matrix, and $f^0$ with the superscript $0$ represents the initial data of $f$. Sometimes we also use
$f_0$ with the subscript $0$ to represent the initial data of $f$.

We denote the Eulerian coordinates by $(x,t)$ with
$x =\zeta(y,t)$ and the Lagrangian coordinates by $(y, t)\in\Omega\times \mathbb{R}_0^+$.
We further assume that, for each fixed $t>0$,
\begin{align}
& \zeta|_{y_3=i}   : \mathbb{R}\to \mathbb{R}\mbox{ is a } C^1(\mathbb{R})\mbox{-diffeomorphism mapping for }i=0,\ 1,\label{20210301715x}\\
&\zeta   :  \overline{\Omega}\to \overline{\Omega} \mbox{ is a } C^1(\overline{\Omega})\mbox{-diffeomorphism mapping},\label{20210301715}
\end{align}
where $\overline{\Omega}:=\mathbb{R}^2\times[0,1]$.

Since $v$ satisfies the divergence-free condition,
and slip boundary condition \eqref{0101}$_5$, we can deduce from \eqref{zeta0inta}--\eqref{zeta0inta0} that
\begin{align}
&\nonumber
J:=\det(\nabla \zeta)=1  \mbox{ in } \Omega,\\
&\nonumber
 \zeta_3 =y_3  \mbox{ on } \partial\Omega .
\end{align}

We introduce the matrix $\mathcal{A}=(\ml{A}_{ij})_{3\times 3}$, which is defined via
\begin{align}\nonumber
\ml{A}^{\mm{T}}=(\nabla\zeta)^{-1}:=
(\partial_j \zeta_i)^{-1}_{3\times 3}.
\end{align}
Then we further define the differential operators $\nabla_{\ml{A}}$, $\mm{div}_{\ml{A}}$ and $\mm{curl}_{\ml{A}}$
as follows: for a given scalar function $f$ and a given vector function $X:=(X_1,X_2, X_3)^{\mm{T}}$,
\begin{align}
&\nabla_{\ml{A}}f:=(\ml{A}_{1k}\partial_kf,
\ml{A}_{2k}\partial_kf, \ml{A}_{3k}\partial_kf)^{\mm{T}},\ \nabla_{\ml{A}}X:=(\nabla_{\ml{A}}X_1,\nabla_{\ml{A}}X_2,\nabla_{\ml{A}}X_3)^{\mm{T}},   \nonumber  \\
&\mm{curl}_{\mathcal{A}}X:=(\ml{A}_{2k}\partial_{k}X_3-\ml{A}_{3k}\partial_{k}X_2,
\ml{A}_{3k}\partial_{k}X_1-\ml{A}_{1k}\partial_{k}X_3,
\ml{A}_{1k}\partial_{k}X_2-\ml{A}_{2k}\partial_{k}X_1
)^{\mm{T}},\nonumber \\
&\mm{div}_{\ml{A}}(X_1,X_2, X_3)^{\mm{T}}:=\ml{A}_{lk}\partial_k X_l, \nonumber
\end{align}
where we have used the Einstein convention of summation over repeated indices, and $\partial_k=\partial_{y_k}$.
In particular, $\mm{curl}X:=\mm{curl}_IX$, where $I$ is a $3\times 3$ identity matrix. In addition, we will denote $(\mm{curl}_{\mathcal{A}} X^1,\ldots, \mm{curl}_{\mathcal{A}}X^n) $ by $\mm{curl}_{\mathcal{A}}(X^1,\ldots X^n) $ for simplicity, where $X^i=(X^i_1,X^i_2,X^i_3)^{\mm{T}}$ is a vector function for $1\leqslant i \leqslant  n$.

Defining the Lagrangian unknowns
\begin{equation*}
(u ,Q, B)(y,t)=(v,P+\lambda|M|^2/2,M)(\zeta(y,t),t) \mbox{ for } (y,t)\in \Omega \times\mathbb{R}^+_0,
\end{equation*}
then in Lagrangian coordinates, the initial-boundary value problem \eqref{0101}  reads as follows:
\begin{equation}\label{01dsaf16asdfasf00}
\begin{cases}
\zeta_t=u,\ \div_\ml{A}u=0 &\mbox{in } \Omega,\\
\rho u_t+\nabla_{\ml{A}}Q+a\rho u=\lambda B\cdot\nabla_{\ml{A}}B &\mbox{in } \Omega,\\
B_t=B\cdot \nabla_{\ml{A}}u  &\mbox{in } \Omega, \\
\div_\ml{A}B=0  &\mbox{in } \Omega, \\
(u,\zeta, B)|_{t=0}=(u^0,\zeta^0, B^0) &\mbox{in } \Omega, \\
(\zeta_3-y_3)=u_3=0  &\mbox{on }\partial\Omega .
\end{cases}
\end{equation}

We can  derive from \eqref{01dsaf16asdfasf00}$_3$ the differential version of magnetic flux conservation \cite{JFJSARMA2019}:
\begin{equation} \nonumber
\ml{A}^{\mm{T}}B=\ml{A}_0^{\mm{T}}B^0,
 \end{equation}
 which yields
\begin{eqnarray}
 \label{0124}  B=\nabla\zeta \ml{A}^{\mm{T}}_0 B^0.
\end{eqnarray}

Let $\bar{M}=(\bar{M}_1,0,0)^{\mm{T}}$ and $\bar{M}_1$
is a nonzero constant.
If we assume
\begin{equation}
\label{201903081437}
\ml{A}^{\mm{T}}_0 B^0=\bar{M}  \mbox{ (i.e.  }B^0=\nabla\zeta^0\bar{M}  = \partial_{\bar{M}}\zeta^0\mbox{)},
 \end{equation}
then \eqref{0124} reduces to
 \begin{eqnarray}
\label{0124xx}   B=\partial_{\bar{M}}\zeta  .
\end{eqnarray}
It should be noted that $B$ given by \eqref{0124xx} automatically satisfies \eqref{01dsaf16asdfasf00}$_3$
and \eqref{01dsaf16asdfasf00}$_4$. Moreover, we see from \eqref{0124xx} that the magnetic tension in Lagrangian coordinates
has the relation
\begin{equation}\nonumber
 \lambda B\cdot \nabla_{\ml{A}} B=  \lambda\bar{M}_1^2\partial_{1}^2\zeta.
 \end{equation}

Let
$$e_1=(1,0,0)^{\mm{T}},\ m= \lambda\bar{M}_1^2 /\rho ,\
q= Q /\rho ,\ \eta=\zeta-y,\ \mathcal{A}=(\nabla \eta+I)^{\mm{T}} \mbox{ and }   \eta^0=\zeta^0-y.$$
Consequently, under the assumption \eqref{201903081437}, the problem \eqref{01dsaf16asdfasf00} is equivalent to
the following initial-boundary value problem:
\begin{equation}\label{01dsaf16asdfasf}
                              \begin{cases}
\eta_t=u &\mbox{ in } \Omega,\\
u_t+\nabla_{\ml{A}}q+a u=m  \partial_1^2\eta  &\mbox{ in } \Omega,\\
\div_\ml{A}u=0  &\mbox{ in } \Omega, \\
(u,\eta)|_{t=0}=(u^0,\eta^0) &\mbox{ in } \Omega, \\
(u_3,\eta_3)=0  &\mbox{on } \partial\Omega ,
\end{cases}
\end{equation}
with the expression of magnetic field
\begin{align}\label{202103190614}
B=\bar{M}_1(\partial_1\eta+e_1).
\end{align}

\subsection{Main results.}\label{subsec:03}

 Before stating our main result, we introduce some notations which will be frequently used throughout this paper.

(1) Simplified basic notations:  $I_a:=(0,a)$, in particular, $I_\infty=\mathbb{R}^+$. $\overline{S}$ denotes the closure of the set $S\subset \mathbb{R}^n$ with $n\geqslant 1$, in particular $\overline{I_T}:=[0,T]$.  $\Omega_a:=\Omega\times I_a$.  $\vec{\mathbf{n}}:=({\mathbf{n}}_1,{\mathbf{n}}_2, {\mathbf{n}}_3)^{\mm{T}}$
denotes the outward unit normal vector on  $\partial\Omega$. $\int:=\int_\Omega=\int_{(0,1)^3}$. $(u)_{\Omega}=\int  u\mm{d}y $.

 $a\lesssim b$ means that $a\leqslant cb$.
If not stated explicitly,  the positive constant $c$ may depend on $a$,  $m $ and $\Omega$, and may vary from one place to other place.
Sometimes we use $c_i$ for $i\geqslant 1$ to replace $c$ in order to emphasize that $c_i$ is fixed value.

 $\alpha$ always denotes the multiindex with respect to the variable $y$,
$|\alpha|=\alpha_1+\alpha_2+\alpha_3$ is called the order of multiindex and  $\partial^{\alpha}:=\partial_{1}^{\alpha_1} \partial_{2}^{\alpha_2} \partial_{3}^{\alpha_3}$.
 ``$\nabla^{i}f\in \mathbb{X}$" represents that $\partial^\alpha f\in \mathbb{X}$ for any multiindex $\alpha$ satisfying $|\alpha|=i$, where $\mathbb{X}$ denotes some set of functions.

(2) Simplified Banach spaces:
\begin{align}
& L^p:=L^p (\Omega)=W^{0,p}(\Omega),\;\;H^{i}:=W^{i,2}(\Omega),
\ L^p_TH^i:=L^p(I_T,H^i),\nonumber \\
&{H}^{l,i}:=\{w\in H^{i}~|~ { \partial_1^{k}w}\in H^{i},\;0\leqslant k\leqslant l\},
\  H^{l,j}_{\mathrm{s}}:=\{w\in {H}^{l,j}~|~w_3|_{\partial\Omega }=0\}, \nonumber \\
&H^{j}_{1}:=\{w \in H^{j} ~|~\det(\nabla w +I )=1\},\
{
H^{l,j}_{1,\mathrm{s}}:= H^{l,j}_{\mathrm{s}}\cap H^{j}_{1}},\nonumber \\
& {{C}}^0_{B,\mm{weak}}(\overline{I_T} ,L^2):= L^\infty_TL^2\cap  C^0(\overline{I_T}, L^2_{\mm{weak}}),  \  \underline{H}^i:=\{w\in H^i~|~(w)_{\Omega}=0\}, \nonumber
\end{align}
where $1\leqslant p\leqslant \infty$,
and  $i$, $l\geqslant 0$ and  $j\geqslant 1$  are integers. In addition,  $\|(f^1,\cdot,f^n)\|_B:=\sqrt{\sum_{1\leqslant i\leqslant n}\|f^i\|_B^2}$, where $f^i$ may be a  scalar function, a vector or a matrix function for $1\leqslant i\leqslant n$.

 (3) Simplified function classes: for any given integer $j\geqslant 1$,
\begin{align}
& H^j_{*}:=\{\xi\in H^j~|~ \xi(y)+y : \overline{\Omega}\to \overline{\Omega}\mbox{ satisfies the properties }\nonumber \\
&\qquad \quad \ \mbox{of  diffeomorphism \eqref{20210301715x} and
\eqref{20210301715} as }\zeta\},\nonumber\\
&\mathfrak{C}^0( \overline{I_T},{H}^{2,3}_{\mathrm{s}}):=\{\eta\in C^0(\overline{I_T} ,{H}^{1,3}_{\mathrm{s}}) ~|~  \nabla^3 \partial_1^2\eta \in   {C}^0_{B,\mm{weak}}(\overline{I_T} ,L^2) \},\nonumber \\
& \underline{\mathfrak{H}}^{2,3}_{1,*,T}:=\{\eta\in \mathfrak{C}^0(\overline{I_T}, H^{2,3}_{\mathrm{s}} )
~|~ \eta(t) \in    H^3_1\cap H^{3}_{*}\mbox{ for each }t\in \overline{I_T}\},\nonumber\\
&    \mathfrak{U}_{T}^{1,3}:=  \{u\in C^0(\overline{I_T}, H^3_{\mathrm{s}}) ~|~\nabla^3 \partial_1 u \in {C}^0_{B,\mm{weak}}( \overline{I_T} ,L^2) , \nonumber   \\
&\qquad \qquad\qquad \qquad\qquad \quad  \
u_t\in C^0(\overline{I_T}, H^2_{\mathrm{s}}),\  u_t\in L^\infty_TH^3  \}\nonumber .
\end{align}

(4) Simplified norms and semi-norms: for integers $i\geqslant 0$ and $j\geqslant 0$,
\begin{align}
\|\cdot \|_{i} :=\|\cdot \|_{H^i},\;\;\;\;
\|\cdot\|_{j,i}:= \|\partial_{1}^{j}\cdot\|_{i},\;\;\;
\ \|\cdot\|_{\underline{j},i}:=\sqrt{\sum_{0\leqslant k\leqslant j}\|\cdot\|_{k,i}^2}.\nonumber
\end{align}

 (5) Simplified functionals:
 \begin{align}&  \mathcal{E}^L:=\|(u,\eta,\partial_1\eta)\|_{3}^2+\|( u_t, \nabla q)\|_2^2,\;\; \mathcal{D}^L:= \|(u,\partial_1\eta)\|_3^2+\|( u_t,\nabla q)\|_2^2\nonumber,\\[2mm]
& \mathcal{E}^{H}:= \|(u, \partial_1\eta )\|_{\underline{1},3}^2
+\|( u_t, \nabla q)\|_2^2,\;\;  \mathcal{D}^{H}:=\|(u,\partial_1u, \partial_1^2\eta)\|_{3}^2+\|( u_t,\nabla q)\|_2^2,\nonumber\\
& I^0_L:= \|(u^0,\eta^0,\partial_1\eta^0)\|_3^2 , \ I^0_H:= \|(u^0,\partial_1 \eta^0)\|_{\underline{1},3}.  \nonumber
               \end{align}
We call $\mathcal{E}^L$, resp. $\mathcal{D}^L$ the lower energy, resp. dissipation functionals, and $\mathcal{E}^{H}$, resp. $\mathcal{D}^H$ the higher energy, resp. dissipation  functionals.

Now we state the existence result of unique global solutions for the initial-boundary value problem \eqref{01dsaf16asdfasf} with   small initial data.

\begin{thm}\label{thm2}
Let
$(u^0,\eta^0)\in H^{1,3}_{\mathrm{s}} \times
 (H^{2,3}_{1,\mathrm{s}}\cap H^3_*)$
satisfy the incompressible condition
$$\mm{div}_{\mathcal{A}^0}u^0=0\mbox{ in }\Omega,$$ where
 $\mathcal{A}^0:=(\nabla \eta^0+I)^{-\mm{T}}$. Then there exist  a sufficiently small constant $\delta >0$ and a positive constant $c$, such that, for any $(u^0,\eta^0)$ satisfying
$$ e^{c \sqrt{I_H^0}}\|(u^0,\eta^0,\partial_1\eta^0)\|_3^2\leqslant\delta^2, $$
 the initial-boundary problem \eqref{01dsaf16asdfasf} admits  a unique global classical solution
$(u,\eta,q)\in \mathfrak{U}_{\mm{s},\infty}^{1,3}\times \underline{\mathfrak{H}}^{2,3,*}_{1,\mm{s},\infty}\times C(\mathbb{R}^+_0 ,\underline{H}^3)$. Moreover the solution enjoys the following stability estimates:
\begin{align}\label{1.200}
&\mathcal{E}^L(t)+\int_0^t\mathcal{D}^L(s)\mm{d}s\lesssim e^{c\sqrt{I^0_H}  }I^0_L,\\
&\label{202101241112}
 {e}^{c_1 t}\mathcal{E}^H(t) + \int_0^t e^{c_1 \tau}\mathcal{D}^H(s)\mm{d}s
\lesssim I^0_H,\\
&e^{c_1 t }\|\eta(t)-\eta^\infty\|_3^2\lesssim {I^0_H},\label{1.200xx}
\end{align}
where $\eta^\infty$ only depends on $y_2$ and $y_3$. In addition, $\nabla_{\mathcal{A}}q\in L^\infty_\infty H^3$.
\end{thm}
\begin{rem}
In this paper, we only consider the impressive horizontal magnetic field $\bar{M}$.
\emph{However it is not clear to the authors that whether there exists a global well-posdness result for non-horizontal magnetic field.}
\end{rem}
\begin{rem}
For the case $\Omega=\mathbb{R}^2\times (0,1)$, we also obtain similar result, where the exponential stability estimates \eqref{202101241112}  and \eqref{1.200xx} should be replaced by  algebraic stability estimates. \emph{We will verify this assertion in a forthcoming paper.}
\end{rem}
\begin{rem}
For each fixed $t\in \mathbb{R}_0^+$, the solution $\eta(y,t)$ in Theorem \ref{thm2} belongs to $  H^3_{*}$. Let $\zeta=\eta+y$, then $\zeta $  satisfies \eqref{20210301715x} and \eqref{20210301715} for each fixed $t\in \mathbb{R}_0^+$. We denote the inverse transformation of $\zeta$ by $\zeta^{-1}$, and then define that
$$
\begin{aligned}
&( v,N,Q)(x,t):=( u(y,t), \bar{M}_1\partial_1\eta(y,t),\rho q(y,t))|_{y=\zeta^{-1}(x,t)},\\
& M:=N+\bar{M}\mbox{ and } P:=Q-\lambda|M|^2/2.
\end{aligned}$$
Consequently  $(  v,M,P)$  is a global solution for \eqref{0101}; moreover, making use of \eqref{1.200}, \eqref{202101241112} and the fifth conclusion in Lemma \ref{20021032019018}, $(v,N,Q)$ enjoys stability estimates, which are similar to \eqref{1.200} and \eqref{202101241112}.
\end{rem}

Next we briefly sketch the proof of Theorem \ref{thm2}, and the details will be presented in Section \ref{sec:global}.

The key proof for the existence of global  small solutions is to derive an  \emph{a priori} energy inequalities \eqref{1.200} and \eqref{202101241112}.
To this purpose,
  let $(\eta,u)$ be a solution to \eqref{01dsaf16asdfasf} and satisfy
\begin{align}
\label{aprpiosesnew}
&\sup\nolimits_{  t\in \overline{I_T}}\|(u,\eta,\partial_1\eta)(t)\|_3\leqslant  {\delta} \in (0,1] \mbox{ for some } T,\\
& \det(\nabla \eta+I)=1,\label{aprpiosasfesnew}
\end{align}
where $\delta$ is sufficiently small.

Under the  \emph{a priori}  assumption  \eqref{aprpiosesnew} of lower-order energy,  we first derive the both energy inequalities of horizontal-type and vorticity-type, see \eqref{202008250856}  and \eqref{202005021632}.
Summing up the two energy inequalities, we can arrive at the lower-order energy inequality  \begin{align}\label{for:0202n0319}
\frac{\mm{d}}{\mm{d}t}\tilde{\mathcal{E}}^L+ \mathcal{D}^L\lesssim (\|u\|_3+\|\nabla q\|_0)\|\eta\|_3^2
\end{align}
for some  energy functional $\tilde{\mathcal{E}}^L$, which is equivalent to $\mathcal{E}^L$.

Obviously we can not directly close the lower-order energy by using \eqref{for:0202n0319}, unless the term ``$\|u\|_3+\|\nabla q\|_0$" decays in time. It should be noted that the term ``$ (\|u\|_3+\|\nabla q\|_0)\|\eta\|_3^2$" (arising from nonlinear terms, see \eqref{20202105112053} and \eqref{2020202105112055}) does not appear in \eqref{for:0202n0319} for the two dimensional cases.  Motivated by the multi-layers energy method introduced in \cite{GYTIAE2}, in which Guo--Tice investigated the global well-posdedness of surface wave problem, we naturally look for another energy inequality, which includes the decay-in-time of $\|u\|_3$.

By careful analysis of energy estimates, we find that, if
 $(u^0,\eta^0)$ additionally satisfies the regularity
 $$\partial_1(u^0,\partial_1\eta^0) \in H^3,$$ the following higher-order energy inequality can be established:
\begin{align}\label{2020122621050319}
\frac{\mm{d}}{\mm{d}t}\tilde{\mathcal{E}}^H + \mathcal{D}^H\leqslant 0,
\end{align}
where $\tilde{\mathcal{E}}^H  $ is equivalent to $\mathcal{E}^H$.
Moreover, thanks to the estimate
$$\|\partial_1\eta\|_3\lesssim\|\partial_1^2\eta\|_3,$$
we have \begin{align}
\nonumber
\mathcal{E}^H \mbox{ is equivalent to }\mathcal{D}^H.
\end{align}
Thus the  higher-order  energy inequality  \eqref{2020122621050319}, together with the above equivalence, immediately implies the exponential stability of higher-order energy \eqref{202101241112},
which particularly presents the exponential decay-in-timely of $ \|u\|_3$ and $\|\nabla q\|_0$. Hence the \emph{a priori }stability estimate \eqref{1.200} further follows  from \eqref{for:0202n0319}.

Thanks to the \emph{a priori} stability estimates \eqref{1.200}--\eqref{202101241112} and the unique local solvability  of the  initial-boundary problem \eqref{01dsaf16asdfasf}, we immediately get the unique global solvability  of the  initial-boundary problem \eqref{01dsaf16asdfasf}.   Finally, \eqref{1.200xx} can be deduced from \eqref{202101241112} by an asymptotic analysis method.

The rest of this paper is organized as follows.
In Sections \ref{sec:global}, we will provide the proof of  Theorem \ref{thm2}. In Section \ref{202102241211}, we will establish the local well-posdedness result for the initial-boundary problem \eqref{01dsaf16asdfasf}.  Finally, in \ref{sec:09}, we will list some well-known mathematical results, which will be used in Sections \ref{sec:global}--\ref{202102241211}.
\section{Proof of Theorem \ref{thm2}}\label{sec:global}
This section is devoted to the proof of Theorem \ref{thm2}.
The key step is to derive the a \emph{priori} stability estimates \eqref{1.200}--\eqref{202101241112} for the initial-boundary value  problem \eqref{01dsaf16asdfasf}.
To this end, let $(u,\eta,q)$ be a solution to the problem \eqref{01dsaf16asdfasf},  and satisfy \eqref{aprpiosesnew} and \eqref{aprpiosasfesnew} with sufficiently small
  $ {\delta}$. Next we  proceed the derivation of \emph{a priori}  estimates.

\subsection{Preliminary estimates}
To begin with, we shall establish some preliminary estimates.
\begin{lem}
\label{201805141072}
Let $0\leqslant i\leqslant 2$ and $1\leqslant j\leqslant 2$.
Under the assumptions \eqref{aprpiosesnew} and \eqref{aprpiosasfesnew}, we have
\begin{enumerate}[(1)]
  \item  the estimates for $\mathcal{A}$ and $\tilde{\mathcal{A}}$:
\begin{align}
&\label{aimdse}
\|\mathcal{A}\|_{C^0(\overline{\Omega})}+ \|\mathcal{A}\|_2 \lesssim  1,\\
&  \label{06041533fwqg}
\|\tilde{\mathcal{A}}\|_{i}\lesssim   \| \eta\|_{i+1}\\
&  \label{06041533fwqgn}
\| \tilde{\mathcal{A}}\|_{j,i}\lesssim   \|  \eta\|_{j,i+1},\\
& \|\ml{A}_t\|_i \lesssim   \|   u\|_{i+1} , \label{06142100x}\\
&  \label{06041safsa533fwqgn}
\| \tilde{\mathcal{A}}_t\|_{1,i}\lesssim   \|  (u,\eta)\|_{1,i+1}.
\end{align}
Here and in what follows $\tilde{\mathcal{A}}:=\mathcal{A}-I$.
\item  the estimates of $\mm{div}{ {u}}$:
\begin{align}
&\label{201808181500} \|\mm{div} u\|_{i}\lesssim \|\eta\|_{3}\|u\|_{i+1},\\
&\label{201808181500n001} \|\mm{div} u\|_{1,i}\lesssim \|(u,\eta )\|_{3}\|(u,\eta )\|_{1,i+1}.
\end{align}
\item  the estimates of $\mm{div} { {\eta}}$:
\begin{align}
&\label{improtian1}
\|\mm{div} { {\eta}}\|_{i} \lesssim \|\eta\|_{3} \|\eta\|_{i+1} ,\\
&\label{improtian1002}
\|\mm{div} { {\eta}}\|_{1,i} \lesssim \|\eta\|_{3} \|\eta\|_{1,i+1} ,\\
&\label{improtian1003}
\|\mm{div} { {\eta}}\|_{2,2} \lesssim \|\eta\|_{3} \|\eta\|_{2,3}+\|\eta\|_{1,3}^2.
\end{align}
\end{enumerate}
\end{lem}
\begin{pf}
 Recalling  the condition \eqref{aprpiosasfesnew} and the definition of $\mathcal{A}$, we have
\begin{equation}
\label{06131422}
\mathcal{A}=(A^*_{ij})_{3\times 3},
\end{equation}
where
$A^{*}_{ij}$ is the algebraic complement minor of $(i,j)$-th entry of matrix $(\partial_j \zeta_i)_{3\times 3}$.  Exploiting \eqref{01dsaf16asdfasf}$_1$, \eqref{aprpiosesnew}, the produce estimate \eqref{fgestims} and Poinc\'are's inequality \eqref{202012241002}, we easily derive \eqref{aimdse}--\eqref{06041safsa533fwqgn} from the expression \eqref{06131422}.
In addition, we easily derive \eqref{201808181500}  and \eqref{201808181500n001} from the incompressible condition $\mm{div}_{\ml{A}}u=0$ by  \eqref{fgestims}.
{ Similarly, we easily get \eqref{improtian1}--\eqref{improtian1003} from the relation  \eqref{aprpiosasfesnew}.}
Interesting readers can refer to \cite{JFJSJMFMOSERT} for the detailed derivation.   \hfill $\Box$
\end{pf}
\begin{lem}\label{lem:07291427nnmmm}
Let $i=0$, $1$, and  the multindex $\alpha$ satisfy $|\alpha|\leqslant2$, $w(t)\in H^3$ for each $t\in \overline{I_T}$
and
\begin{align}
\label{2021011281929}
{W}^{i,\alpha}:= \partial_{1}^{i}\partial^{\alpha}\left(\mm{curl}_{\mathcal{A}_{t}}u
-m \mm{curl}_{\partial_1\mathcal{A}} \partial_1\eta\right).
\end{align}
Under the assumptions \eqref{aprpiosesnew} and \eqref{aprpiosasfesnew}, we have
\begin{align}
&\label{202101231432}
\left\|\partial^{\alpha}\left(\operatorname{curl}_{\mathcal{A}} w-\operatorname{curl}w\right)\right\|_0
\lesssim\|\eta\|_3\|w\|_3,\\
&\label{202101231245}
\left\|\partial^{\alpha} \left(\operatorname{curl}_{\mathcal{A}} w-\operatorname{curl}w\right)\right\|_{1,0}\lesssim
 \|(\eta,w)\|_3 \|(\eta,w)\|_{1,3},\\
&\label{201907291200000}
\left\|W^{0,\alpha}\right\|_0\lesssim\|(u,\partial_1\eta)\|_3^2 ,\\
&\label{201907291200000nm}
\left\|W^{1,\alpha}\right\|_0\lesssim\|u\|_3\|(u,\eta)\|_{1,3}
+\| \eta\|_{1,3}\| \eta\|_{2,3} .
\end{align}
\end{lem}
\begin{pf}
By  \eqref{06041533fwqg}, \eqref{06041533fwqgn} and the product estimate \eqref{fgestims}, we easily get \eqref{202101231432} and \eqref{202101231245}.
Similarly, making use of \eqref{06041533fwqgn}--\eqref{06041safsa533fwqgn} and \eqref{fgestims}, we easily derive \eqref{201907291200000}--\eqref{201907291200000nm}  from the definition of $W^{i,\alpha}$.
The proof is complete.
\hfill$\Box$
\end{pf}
\subsection{Horizontal spatial estimates}
Now we derive the estimates of horizontally spatial derivatives of $\left(u,\eta\right)$.
\begin{lem}\label{lem:08241445}
Let $0\leqslant k\leqslant1$. Under the assumptions \eqref{aprpiosesnew} and \eqref{aprpiosasfesnew},
we have
\begin{align}
&\label{202008241448}
\frac{\mm{d}}{\mm{d}t}\left(\|\partial_1^{k} u\|^2_0
+m \|\partial_1^{k+1} \eta\|_{0}^2 \right)+ a\|\partial_1^{k} u\|_{0}^2
\lesssim
\begin{cases}
0
&\hbox{for }k=0; \\
\sqrt{\mathcal{E}^L}\mathcal{D}^H
&\hbox{for }k=1,
\end{cases}
\\
&
\frac{\mm{d}}{\mm{d}t}\left(\int\partial_1^{k}\eta\cdot\partial_1^{k} u\mm{d}y
+\frac{a}{2}\|\partial_1^{k} \eta\|_{0}^2\right)
+m \| \partial_1^{k+1} \eta\|_{0}^2
\nonumber \\
&\lesssim\|\partial_1^{k}u\|_{0}^2
+
\begin{cases}
{ \|\eta\|_2^2\|\nabla q\|_0}
&\hbox{for }k=0; \\
\sqrt{\mathcal{E}^L}\mathcal{D}^H
&\hbox{for } k=1\label{202008241446}.
\end{cases}
\end{align}
\end{lem}
\begin{pf}
Applying $\partial_1^{k}$ to \eqref{01dsaf16asdfasf} yields
\begin{equation}\label{01dsaf16asdfasf0321}
                              \begin{cases}
\partial_1^{k}\eta_t=\partial_1^{k}u &\mbox{in } \Omega,\\
\partial_1^{k}u_t+\partial_1^{k}\nabla_{\ml{A}}q+a \partial_1^{k}u=m  \partial_1^2\partial_1^{k}\eta  &\mbox{in } \Omega,\\
\partial_1^{k}\div_\ml{A}u=0  &\mbox{in } \Omega, \\
\partial_1^{k}(u_3, \eta_3)=0  &\mbox{on } \partial\Omega ,
\end{cases}
\end{equation}
Multiplying \eqref{01dsaf16asdfasf0321}$_2$ by $\partial_1^{k}u$ and $\partial_1^{k}\eta$ in $L^2$, resp., and then using
the integration by parts and \eqref{01dsaf16asdfasf}$_1$,
we have
\begin{align}
&\label{202008241510n}
\frac{1}{2}\frac{\mm{d}}{\mm{d}t}\left(\|\partial_1^{k} u\|_0^2
+m \|\partial_1\partial_1^{k}\eta\|_0^2\right)+a\|\partial_1^{k} u\|_{0}^2
=-\int \partial_1^{k}\nabla_{\ml{A}}q \cdot\partial_1^{k}u\mm{d}y
\end{align} and
\begin{align}
&\label{202008241510}
\frac{\mm{d}}{\mm{d}t}\left(\int\partial_1^{k}\eta\cdot\partial_1^{k} u\mm{d}y
+\frac{a}{2}\|\partial_1^{k} \eta\|_{0}^2\right)
+m \|\partial_1\partial_1^{k}\eta\|_0^2
=\|\partial_1^{k}u\|_{0}^2
-\int \partial_1^{k}\nabla_{\ml{A}}q\cdot\partial_1^{k} \eta\mm{d}y.
\end{align}
Next we estimate for the integrals involving the pressure in \eqref{202008241510n} and \eqref{202008241510} by two cases.

(1) \emph{Case $k=0$.}

Noting that $\eta_3|_{\partial\Omega}=0$,
 thus
\begin{align}
\mathcal{A}_{13}|_{\partial\Omega}=\mathcal{A}_{23}|_{\partial\Omega}=0.
\label{202010230522238}
\end{align}
Exploiting the integration by parts, \eqref{01dsaf16asdfasf}$_3$  and the  boundary-value conditions of $u_3|_{\partial\Omega}=0$ and \eqref{202010230522238}, we have
\begin{align}\nonumber
-\int\nabla_{\ml{A}}q\cdot u\mm{d}y=\int q\mm{div}_{\ml{A}}u\mm{d}y=0.
\end{align}
Putting the above identity into \eqref{202008241510n} with $k=0$ yields \eqref{202008241448} with $k=0$.

Now we turn to estimating for the last integral with $k=0$ in \eqref{202008241510}.
Since $\det(\nabla \eta+I)=1$,  we have
\begin{align}
\label{202103051850}
\mm{div}\eta=\mm{div}\Psi(\eta),
\end{align} where
\begin{align}\nonumber
  \Psi(\eta):=-\left(\begin{array}{c}
          \eta_1(\partial_2\eta_2+\partial_3\eta_3 )-
           \eta_1(\partial_2\eta_3 \partial_3\eta_2
- \partial_2\eta_2\partial_3\eta_3) \\
         \eta_2\partial_3\eta_3-\eta_1\partial_1\eta_2-
        \eta_1( \partial_1\eta_2 \partial_3\eta_3
-\partial_1\eta_3\partial_3\eta_2)\\
         -\eta_1\partial_1\eta_3
-\eta_2\partial_2\eta_3-\eta_1(\partial_1\eta_3 \partial_2\eta_2
-\partial_1\eta_2\partial_2\eta_3)
        \end{array}\right).
\end{align}
It is easy to estimate that
\begin{equation}\label{201910040902n00}
\|\Psi(\eta)\|_0\lesssim  \|\eta\|_2^2.
\end{equation}

Thanks to \eqref{202103051850}, \eqref{201910040902n00} and the boundary-value conditions $$\eta_3|_{\partial\Omega}=0,\ \Psi_3(\eta)|_{\partial\Omega}=0,$$
we have
\begin{align}
\int q\mm{div} \eta\mm{d}y
\leqslant\|\Psi(\eta)\|_0\|\nabla q\|_0
\lesssim \|\eta\|_2^2\|\nabla q\|_0. \nonumber
\end{align}
In addition, we can estimate that
\begin{equation}\nonumber
\left|\int \nabla_{\tilde{\ml{A}}}q\cdot\eta\mm{d}y\right|
\lesssim \|\eta\|_2^2\|\nabla q\|_0.
\end{equation}
Exploiting  the above two estimates and integration by parts, we get
\begin{align}
-\int\nabla_{\ml{A}}q\cdot\eta\mm{d}y
 =\int ( q\mm{div}\eta- \nabla_{\tilde{\ml{A}}}q\cdot\eta)\mm{d}y\lesssim \|\eta\|_2^2\|\nabla q\|_0. \label{20202105112053} \end{align}
 Putting the above estimate into \eqref{202008241510} yields  \eqref{202008241446} with $k=0$.

(2) \emph{Case $k=1$.}

Thanks to \eqref{06041533fwqg}, \eqref{201808181500n001} and \eqref{202012241002},  we can estimate that
\begin{align}
&-\int\partial_1\nabla_{\ml{A}}q\cdot\partial_1u\mm{d}y
 =\int(\partial_1q\partial_1 \mm{div}u
+ \nabla_{\tilde{\ml{A}}}q\cdot\partial_1^2u)\mm{d}y\nonumber\\
&\lesssim
\|(u, \eta)\|_3\|(u,\eta)\|_{1,1}\|\partial_1 q\|_0
+\|\partial_1^2u\|_0\| \eta\|_3 \|\nabla q\|_0\lesssim
\sqrt{\mathcal{E}^L}\mathcal{D}^H. \nonumber
\end{align}
Putting the above estimate into \eqref{202008241510n} yields \eqref{202008241448} with $k=1$.

Thanks to \eqref{06041533fwqg}, \eqref{improtian1002}   and \eqref{202012241002}, we have
\begin{align}
&-\int\partial_1\nabla_{\ml{A}}q \cdot\partial_1 \eta\mm{d}y
=\int( \nabla_{\tilde{\ml{A}}}q \cdot\partial_1^2 \eta
+ \partial_1q \mm{div}\partial_1 \eta)\mm{d}y\nonumber \\
&\lesssim  \| \eta\|_{2,0} \|\eta\|_3\|\nabla q\|_0
+\| \eta\|_{1,1}\|\eta\|_3\|  q\|_{1}
\lesssim \sqrt{\mathcal{E}^L}\mathcal{D}^H. \nonumber
\end{align}
Plugging the above estimate into \eqref{202008241510} yields \eqref{202008241446} with $k=1$.
The proof is complete.
\hfill $\Box$
\end{pf}

Lemma \ref{lem:08241445} immediately implies the following two horizontal-type energy  inequalities.
\begin{pro}\label{pro0902}
Under the assumptions \eqref{aprpiosesnew} and \eqref{aprpiosasfesnew}, there exist two energy functionals $\mathcal{E}^H_{h}$ and $\mathcal{E}^L_{h}$ of $(\eta,u)$
 such that
\begin{align}
&\label{202008250856n0}
\frac{\mm{d}}{\mm{d}t}
\mathcal{E}^L_{h}
+  \|(u,\partial_1\eta)\|_{0}^2
\lesssim \|\eta\|_2^2 \|\nabla q\|_0 ,\\
&\label{202008250856}
\frac{\mm{d}}{\mm{d}t}
\mathcal{E}^H_{h}
+ \|(u,\partial_1u, \partial_1^2\eta)\|_{0}^2 \lesssim\sqrt{\mathcal{E}^L}\mathcal{D}^H,
\end{align}
where, for each $t>0$,
$$
 \mathcal{E}^L_{h}\mbox{ and } \mathcal{E}^H_{h}
\mbox{ are equivalent to }
\|(u,\eta,\partial_1\eta )\|_{0}^2\mbox{ and }\|(u,\partial_1\eta)\|_{\underline{1},0}^2 \mbox{, resp.}.
$$
\end{pro}

\subsection{Vorticity estimates}
Let the multindex $\alpha $ satisfy $|\alpha|\leqslant 2$.
Applying $\partial^{\alpha  } \partial_1^{i}\mm{curl}_{\mathcal{A}} $
to \eqref{01dsaf16asdfasf}$_2$ yields
\begin{align}\label{202005021542000}
\partial_t\partial^{\alpha}{ \partial_1^{i}}\mm{curl}_{\mathcal{A}}u
+a\partial^{\alpha}{ \partial_1^{i}}\mm{curl}_{\mathcal{A}}u
=m \partial_1(\partial^{\alpha}{ \partial_1^{i}}\mm{curl}_{\mathcal{A}}\partial_1\eta)
+W^{i,\alpha},
\end{align}
where $W^{i,\alpha}$ is defined by \eqref{2021011281929}.

Now we derive the following vorticity estimates for $u$ and $ \eta$.
\begin{lem}\label{2019100216355nnn00}
Let $|\alpha|\leqslant 2$. Under the assumptions \eqref{aprpiosesnew} and \eqref{aprpiosasfesnew}, we have
\begin{align}
&\label{202005021610nnmm}
\frac{1}{2}\frac{\mm{d}}{\mm{d}t}
\int\left(|\partial^{\alpha}{ \partial_1^i}\mm{curl}_{\mathcal{A}}u |^2+m |\partial^{\alpha}{ \partial_1^i}\mm{curl}_{\mathcal{A}}\partial_1\eta|^2\right)\mm{d}y
 +a\|\partial^{\alpha}{ \partial_1^i}\mm{curl}_{\mathcal{A}}u\|^2_{0}\nonumber\\
&\lesssim{
\begin{cases}
 \|u\|_3 \|(u,{ \partial_1\eta})\|_3^2
&\hbox{for }   i=0; \\
  \| \eta\|_{2,3}(\|u\|_3\| \eta\|_{2,3}+ \| ( u,\eta)\|_{1,3}\|\eta\|_{1,3}
)
  +  \| u\|_3  \| (u,\eta)\|_{1,3}^2
&\hbox{for }  i=1,
\end{cases}}\\
&\label{202005021600nnmm}
\frac{\mm{d}}{\mm{d}t}
\int\left(\partial^{\alpha}{ \partial_1^i}\mm{curl}_{\mathcal{A}}\eta \cdot\partial^{\alpha}{ \partial_1^i}\mm{curl}_{\mathcal{A}} u
+\frac{a}{2}|\partial^{\alpha}{ \partial_1^i}\mm{curl}_{\mathcal{A}} \eta |^2\right)\mm{d}y
+{m }\|\partial^{\alpha}{ \partial_1^i}\mm{curl}_{\mathcal{A}}\partial_1\eta\|^2_{0}\nonumber\\
&\lesssim\|\partial^{\alpha}{ \partial_1^i}\mm{curl}_{\mathcal{A}}u\|^2_{0}\nonumber \\
&\quad +
\begin{cases}
\|\eta\|_3 \|(u,\partial_1\eta)\|_3^2
+{\|u\|_3  \|\eta\|_3^2}
&\hbox{for }   i=0; \\
  \|(u,\eta)\|_3   \| (u,\eta)\|_{1,3}^2
+ \| \eta\|_{2,3}(\|\eta\|_3\| \eta\|_{2,3} +\|\eta\|_{1,3}^2)
&\hbox{for }  i=1.
\end{cases}
\end{align}
\end{lem}
\begin{pf}
(1) Multiplying \eqref{202005021542000} by $\partial^{\alpha}{ \partial_1^i}\mm{curl}_{\ml{A}}u$ in $L^2$ and then using the integration by parts, we get
\begin{align}
&\frac{1}{2}\frac{\mm{d}}{\mm{d}t}
\int\left(|\partial^{\alpha}{ \partial_1^i}\mm{curl}_{\mathcal{A}}u |^2+m |\partial^{\alpha}{ \partial_1^i}\mm{curl}_{\mathcal{A}}\partial_1\eta|^2\right)\mm{d}y
 +a\|\partial^{\alpha}{ \partial_1^i}\mm{curl}_{\mathcal{A}}u\|^2_{0}\nonumber \\
&= m \int\partial^{\alpha}{ \partial_1^i}\mm{curl}_{\mathcal{A}} \partial_1\eta\cdot\partial^{\alpha}{ \partial_1^i}\mm{curl}_{\mathcal{A}_t} \partial_1\eta\mm{d}y-m \int\partial^{\alpha}{ \partial_1^i}\mm{curl}_{\mathcal{A}}\partial_1\eta
\cdot\partial^{\alpha}{ \partial_1^i}\mm{curl}_{\partial_1\mathcal{A}} u\mm{d}y
\nonumber
\\
&\quad +\int{W^{i,\alpha}}
\cdot\partial^{\alpha}{ \partial_1^i}\mm{curl}_{\mathcal{A}} u \mm{d}y=:\sum_{j=1}^{3} I_j^{i,\alpha}.
\label{202005031742nnmm} \end{align}

Making use of
\eqref{aimdse}, \eqref{06041533fwqgn}, \eqref{06142100x}, and \eqref{201907291200000},
the three integrals $I_1^{0,\alpha} $--$I_3^{0,\alpha} $ can be estimated as follows:
\begin{align}
I_1^{0,\alpha}
\lesssim&\|\mm{curl}_{\mathcal{A}}\partial_1\eta\|_2 \|\mm{curl}_{\mathcal{A}_t} \partial_1\eta\|_2
\lesssim   \|u\|_3\| \eta\|_{1,3}^2 ,\label{202105021842}\\
 I_2^{0,\alpha} \lesssim& \|\mm{curl}_{\mathcal{A}}\partial_1\eta\|_2\| \mm{curl}_{\partial_1\mathcal{A}} u \|_2\lesssim  \|u\|_3\| \eta\|_{1,3}^2,  \\
 I_3^{0,\alpha}\lesssim&\|W^{0,\alpha}\|_0\|\mm{curl}_{ \mathcal{A}} u\|_2
\lesssim \|u\|_3 \|(u,\partial_1\eta)\|_3^2.
\end{align}

Exploiting \eqref{aprpiosesnew},
 \eqref{aimdse}, \eqref{06041533fwqgn}--\eqref{06041safsa533fwqgn}, \eqref{201907291200000nm} and \eqref{202012241002},
the three integrals $I_1^{1,\alpha}$--$I_3^{1,\alpha}$ can be estimated as follows:
\begin{align}
I_1^{1,\alpha}&\lesssim \| \eta\|_{2,3} (\|u\|_3\|\eta\|_{2,3}+\|(u,\eta)\|_{1,3}
\|\eta\|_{1,3}) , \\
I_2^{1,\alpha}
& \lesssim \|\eta\|_{2,3}(\|u\|_3\|\eta\|_{2,3}+\| u\|_{1,3}\|\eta\|_{1,3} ) ,\\
I_3 ^{1,\alpha}
&\lesssim  \|(u,\eta)\|_{1,3}(\|u\|_3 \|(u,\eta)\|_{1,3}+\| \eta\|_{1,3}\|\eta\|_{2,3}) .\label{202105021842xx}
\end{align}
Thanks to the above six estimates and  \eqref{202012241002}, we immediately get \eqref{202005021610nnmm} from \eqref{202005031742nnmm}.

(2) Multiplying \eqref{202005021542000}
by $\partial^{\alpha}{ \partial_1^i}\mm{curl}_{\ml{A}}\eta$ in $L^2$, and then using
the integration by parts, we obtain
\begin{align}
&\frac{\mm{d}}{\mm{d}t}
\int\left(\partial^{\alpha}{ \partial_1^i}\mm{curl}_{\mathcal{A}}\eta \cdot\partial^{\alpha}{ \partial_1^i}\mm{curl}_{\mathcal{A}} u
+\frac{a}{2}|\partial^{\alpha}{ \partial_1^i}\mm{curl}_{\mathcal{A}} \eta |^2\right)\mm{d}y
 +{m }\|\partial^{\alpha}{ \partial_1^i}\mm{curl}_{\mathcal{A}}\partial_1\eta\|^2_{0}\nonumber \\
&
=\|\partial^{\alpha}{ \partial_1^i}\mm{curl}_{\mathcal{A}}u\|^2_{0}+
\int \partial^{\alpha}{ \partial_1^i}\mm{curl}_{\mathcal{A}_t}\eta \cdot\partial^{\alpha}{ \partial_1^i}\mm{curl}_{\mathcal{A}} u\mm{d}y
+{a}\int\partial^{\alpha}{ \partial_1^i}\mm{curl}_{\mathcal{A}_t} \eta
\cdot\partial^{\alpha}{ \partial_1^i}\mm{curl}_{\mathcal{A}} \eta  \mm{d}y\nonumber \\
&\quad -m \int (\partial^{\alpha}{ \partial_1^i}\mm{curl}_{\mathcal{A}} \partial_1\eta)
\cdot\partial^{\alpha}{ \partial_1^i}\mm{curl}_{\partial_1\mathcal{A}} \eta \mm{d}y+\int
{ W^{i,\alpha}}\cdot\partial^{\alpha}{ \partial_1^i}\mm{curl}_{\mathcal{A}} \eta\mm{d}y
\nonumber \\
&=:\|\partial^{\alpha}{ \partial_1^i}\mm{curl}_{\mathcal{A}}u\|^2_{0}{ + \sum_{j=4}^7 I_j^{i,\alpha}}.
\label{202005031720nnmm}
\end{align}

Similarly to \eqref{202105021842}--\eqref{202105021842xx},
the   integrals $I_4^{i,\alpha} $--$I_7^{i,\alpha} $  can be estimated as follows:
\begin{align}
& I_4^{0,\alpha} \lesssim \|u\|_3^2\|\eta\|_3,\
 I_5^{0,\alpha}
\lesssim \|u\|_3\|\eta\|_3^2 , \label{2020202105112055} \\
&  I_6^{0,\alpha} \lesssim
 \|\eta\|_3\| \eta\|_{1,3}^2 ,\ I_7^{0,\alpha}
 \lesssim  \|\eta\|_3 \|(u,\partial_1\eta)\|_3^2, \nonumber \\
 &I_4^{1,\alpha},\ I_5^{1,\alpha}
\lesssim   \|(u,\eta)\|_3  \|(u,\eta)\|_{1,3}^2  , \nonumber
\\
&I_6^{1,\alpha} \lesssim   \| \eta\|_{2,3} ( \|\eta\|_3\|\eta\|_{2,3} +\|\eta\|_{1,3}^2)
 , \nonumber  \\
&I_7^{1,\alpha} \lesssim   \|\eta\|_{1,3} (\|u\|_3\|(u,\eta)\|_{1,3}
+\|\eta\|_{1,3}\|\eta\|_{2,3}). \nonumber
\end{align}
Thanks to the above eight estimates,
we immediately
get \eqref{202005021600nnmm} from \eqref{202005031720nnmm}.
This completes the proof of Lemma \ref{2019100216355nnn00}.
\hfill$\Box$
\end{pf}

Next we further establish the vorticity-type energy inequalities.
\begin{pro}\label{pro1632nn1}
Under the assumptions \eqref{aprpiosesnew} and \eqref{aprpiosasfesnew}, there exist two energy functionals $\mathcal{E}^H_{\mm{cul}}$ and $\mathcal{E}^L_{\mm{cul}}$ of $(\eta,u)$
 such that
\begin{align}
\label{202005021632nn1}
&\frac{\mm{d}}{\mm{d}t}
\mathcal{E}^L_{\mm{cul}}
+  \|\mm{curl}( u,\partial_1 \eta)\|_2^2
\lesssim \|(u,\eta)\|_3\|(u,\partial_1\eta)\|_3^2
+ \|u\|_3\|\eta\|_3^2
\end{align}
and
\begin{align}\label{202005021632}
&\frac{\mm{d}}{\mm{d}t}
\mathcal{E}^H_{\mm{cul}}
+ \|\mm{curl}(u, \partial_1 u,\partial_1^2  \eta)\|_{2 }^2
\lesssim \|(u,\eta,\partial_1\eta)\|_3 \|(u, \partial_1\eta )\|_{\underline{1},3}^2,
\end{align}
where, for each $t\in \overline{I_T}$,
\begin{align}
& \mathcal{E}^L_{\mm{cul}}
\mbox{ is equivalent to }
 \|\mm{curl}_{\mathcal{A}}(u,\eta, \partial_1\eta )\|_{2}^2,\label{2022100520202157} \\
&\mathcal{E}^H_{\mm{cul}}
\mbox{ is equivalent to }
{  \|  \mm{curl}_{\mathcal{A}}(u,\partial_1\eta)
\|_{\underline{1},2}^2+\|\mm{curl}_{\mathcal{A}}\eta\|_{1,2}^2}.\label{2022100520202157x}
\end{align}
\end{pro}
\begin{pf}
We easily derive from \eqref{202005021610nnmm} and \eqref{202005021600nnmm} with $i=0$ that
\begin{align}
\label{202101242005}
&\frac{\mm{d}}{\mm{d}t}
\mathcal{E}^L_{\mm{cul}}
+  \| \mm{curl}_{\mathcal{A}}(u,\partial_1\eta)\|_2^2
\lesssim
\|(u,\eta)\|_3\|(u,\partial_1\eta)\|_3^2
+ \|u\|_3\|\eta\|_3^2,
\end{align}
where  $\mathcal{E}^L_{\mm{cul}}$ satisfies  \eqref{2022100520202157}.
Employing \eqref{202101231432} with $w=u$ and $ \partial_1\eta$, we further get \eqref{202005021632nn1}
from \eqref{202101242005}.

 We derive from \eqref{202005021610nnmm} and \eqref{202005021600nnmm} with $i=1$ that
\begin{align}
\label{202005021632nx2}
&\frac{\mm{d}}{\mm{d}t}
\mathcal{E}^{H}_{\mm{cul}}
+\|\mm{curl}_{\mathcal{A}}u\|^2_{2} +  \|\mm{curl}_{\mathcal{A}}\left(u
, \partial_1\eta\right)\|^2_{1,2}
\lesssim
\|(u,\eta,\partial_1\eta)\|_3(\|u\|_3^2+\| (u, \eta,\partial_1 \eta)\|_{1,3}^2),
\end{align}
where  $\mathcal{E}^{H}_{\mm{cul}}$ satisfies \eqref{2022100520202157x}.
Consequently, making use of \eqref{202101231432} and \eqref{202101231245}, we derive \eqref{202005021632nn1} and \eqref{202005021632}
from \eqref{202101242005} and \eqref{202005021632nx2}.
This completes the proof of Proposition \ref{pro1632nn1}.
\hfill$\Box$
\end{pf}

\subsection{Equivalent estimates }\label{subsec:equivalent}

This section is devoted to establishing the equivalent estimates for
 $\mathcal{E}^L$, $\mathcal{D}^L$, $\mathcal{E}^H$ and $\mathcal{D}^H$.
\begin{lem}\label{201612132242nx}
Under assumptions \eqref{aprpiosesnew} and \eqref{aprpiosasfesnew}, we have
\begin{align}
&\label{201903161543}
\mathcal{E}^L\mbox{ is equivalent to }\|(u,\eta,\partial_1 \eta )\|_3^2,\\
&\label{202103192137}
\mathcal{D}^L\mbox{ is equivalent to }\|(u,\partial_1 \eta )\|_3^2,\\
&\mathcal{E}^H,\ \mathcal{D}^H\mbox{ and }\|(u,\partial_1u,\partial_1^2\eta)\|_{3}^2\mbox{ are equivalent}.
\label{2019031615431}
\end{align}
\end{lem}
\begin{pf}To obtain \eqref{201903161543}--\eqref{2019031615431}, the key step is to establish the estimates of $\nabla q$ and $u_t$.
By  \eqref{01dsaf16asdfasf}$_2$, we have
\begin{equation}
\begin{cases}
-\Delta q= f^1  &\mbox{in } \Omega,\\
\partial_3 q =f^2  \vec{\mathbf{n}}_3 &\mbox{on } \partial\Omega,
\end{cases}
\label{2022010502022007}
\end{equation}
where $f^1:=\mm{div}\left(u_t+a u-m \partial_1^2\eta+\nabla_{\tilde{\ml{A}}}q\right)$ and
$f^2:=- {\tilde{\ml{A}}}_{3j}\partial_j q$.

Note that
\begin{align*}
\int f^1\mm{d}y
+\int_{\partial\Omega} f^2  { \vec{\mathbf{n}}_3}\mathbf{d}y_{\mm{h}}=0,
\end{align*}
  where $\vec{\mathbf{n}}_3$ is the third component of the outward unit normal vector $\vec{\mathbf{n}} $ on  $\partial\Omega$, thus applying  the elliptic estimate \eqref{neumaasdfann1n} to \eqref{2022010502022007}, and then using \eqref{01dsaf16asdfasf}$_3$, we have
\begin{align}
\|\nabla q\|_2&\lesssim\|f^1\|_1+\|f^2\|_2 \nonumber \\
&\lesssim \|(\mm{div}_{\ml{A}_t}u,\mm{div}_{\tilde{\ml{A}}}u_t,\mm{div}_{\tilde{\ml{A}}}u)\|_1
+\|(\partial_1^2\eta, \nabla_{\tilde{\ml{A}}}q)\|_2 \nonumber \\
&\lesssim\| \eta\|_{2,2}
+ \|(u,\eta)\|_3\|(u,\nabla q,u_t)\|_2. \label{202008241115}
\end{align}

In addition, by \eqref{01dsaf16asdfasf}$_2$, we have
\begin{align}\nonumber
\|u_t\|_2
\lesssim\|(u,\partial_1^2\eta,\nabla q)\|_2,
\end{align}
which, together with \eqref{202008241115}, yields, for sufficiently small $\delta$,
\begin{align}\label{202008241115nmq}
\|(u_t,\nabla q)\|_2
\lesssim\|(u,\partial_1^2\eta )\|_2.
\end{align}
Consequently, we conclude the desired results by   \eqref{202008241115nmq}
and  \eqref{202012241002}.
The proof is complete.
\hfill$\Box$
\end{pf}

\subsection{Stability estimates}

Now we are in the position to  derive  \emph{a priori} stability estimates \eqref{1.200} and \eqref{202101241112}.

 Making use of \eqref{201808181500}--\eqref{improtian1003},   and Hodge-type elliptic estimate \eqref{08171537}, we can deduce  that, for  $0\leqslant i\leqslant 1$ and $0\leqslant j\leqslant 2$,
\begin{align}
\| u \|_{i,3}\lesssim &\|  u \|_{i,0}
+\| \mm{curl} u  \|_{i,2} +\begin{cases}
\|\eta\|_{3} \|u\|_{3}&\mbox{for }  i=0 ,
\\
\|(u,\eta)\|_3\|( u,\eta)\|_{1,3} &\mbox{for }i=1   \end{cases} \label{202008240900} \end{align}
and
\begin{align}
\|  \eta \|_{j,3}\lesssim &\|  \eta \|_{j,0}
+\| \mm{curl}  \eta \|_{j,2}  +  \begin{cases}
\|\eta\|_{3}^2&\mbox{for }  j=0 ;
\\
\|\eta\|_{3} \|\eta\|_{1,3} &\mbox{for }  j=1; \\
\|\eta\|_{3} \|\eta\|_{2,3}+ \|\eta\|_{1,3}^2 &\mbox{for }  j=2. \end{cases}\label{20200824sfas0900} \end{align}

By \eqref{202008250856}, \eqref{202005021632} and \eqref{202012241002}, we have
\begin{align}\label{201910050940}
\frac{\mm{d}}{\mm{d}t}\tilde{\mathcal{E}}^H+c\tilde{\mathcal{D}}^H\lesssim \sqrt{\mathcal{E}^L}\mathcal{D}^H,
\end{align}
where
$$\begin{aligned}
&\tilde{\mathcal{E}}^H:=
\mathcal{E}^H_h +\mathcal{E}^H_{\mm{cul}} \mbox{ and }
 \tilde{\mathcal{D}}^H:=
\|(u,\partial_1u,\partial_1^2\eta)\|_{0}^2
+\|\mm{curl}(u,\partial_1 u,\partial_1^2 \eta )\|_2^2.
\end{aligned}$$
In addition, it follows from \eqref{202101231432}, \eqref{202101231245}, \eqref{2019031615431},
\eqref{202008240900} and \eqref{20200824sfas0900}
that for sufficiently small ${\delta}$,
\begin{align}
&\label{201908081232}
\tilde{\mathcal{E}}^H,\ {\mathcal{E}^H},\ \tilde{\mathcal{D}}^H,\ {\mathcal{D}^H}\mbox{ and }\|(u,\partial_1u, \partial_1^2\eta)\|^2_{3}
\mbox{ are equivalent}.
\end{align}
Thus  we deduce from \eqref{201910050940}--\eqref{201908081232} that,
\emph{for any sufficiently small $\delta$,}
\begin{align}
\nonumber
\frac{\mm{d}}{\mm{d}t}\tilde{\mathcal{E}}^H+\mathcal{D}^H\leqslant 0.
\end{align}

Exploiting \eqref{201908081232}, we further deduce from the above inequality that  there exists a constant  $c_1>0$ such that
\begin{align}    \label{estemalas}
 & {e}^{c_1 t}\mathcal{E}^H(t) + \int_0^t e^{c_1 s}\mathcal{D}^H(s)\mm{d}{s}
\lesssim I_H^0.
\end{align}

In addition, by \eqref{202008250856n0} and \eqref{202005021632nn1}, we have
\begin{align}\label{201910050940nnm}
\frac{\mm{d}}{\mm{d}t}\tilde{\mathcal{E}}^L+c\tilde{\mathcal{D}}^L\lesssim  \|(u,\eta)\|_3 \|(u, \partial_1\eta )\|_{3}^2 + (\|u\|_3+\|\nabla q\|_0)\|\eta\|_3^2,
\end{align}
where
$$\begin{aligned}
 \tilde{\mathcal{E}}^L:=
\mathcal{E}^L_h+\mathcal{E}^L_{\mm{cul}}\mbox{ and }\tilde{\mathcal{D}}^L:=
 \|( u,\partial_1\eta )\|_{0}^2
 +\|\mm{curl}(u,\partial_1\eta)\|_{2}^2.
\end{aligned}$$
In addition, thanks to \eqref{202101231432},
  \eqref{201903161543}, \eqref{202103192137}, \eqref{202008240900} and \eqref{20200824sfas0900}, we have, for sufficiently small $\delta$,
\begin{align}
&\label{2019080sdaf81232nnm}
\tilde{\mathcal{E}}^L,\ {\mathcal{E}}^L\mbox{ and }\|( u,\eta,\partial_1\eta)\|_{3}^2
\mbox{ are equivalent}
\end{align}
 and
\begin{align}
 \label{201908081232nnm}
\tilde{\mathcal{D}}^L\mbox{ is equivalent to }{\mathcal{D}^L}.
\end{align}

Thanks to \eqref{201908081232nnm}, we get from \eqref{201910050940nnm} that,  \emph{for sufficiently small $\delta$},
\begin{align}\nonumber
\frac{\mm{d}}{\mm{d}t}\tilde{\mathcal{E}}^L+c\tilde{\mathcal{D}}^L\lesssim   (\|u\|_3+\|\nabla q\|_0)\|\eta\|_3^2.
\end{align}
Exploiting \eqref{estemalas}, \eqref{2019080sdaf81232nnm}   and Gronwall's lemma, we deduce from the above estimate  that
\begin{align} \label{202101242212}
 \mathcal{E}^L(t) +\int_0^t \mathcal{D}^L(s)\mm{d}{s} \lesssim
e^{c  \sqrt{I_H^0}}I_L^0.
\end{align}
In particular, by \eqref{estemalas} and \eqref{202101242212}, there exist constants $c_2\geqslant 1$ and $\delta_1\in (0,1]$  such that, for any $\delta\leqslant \delta_1$,
\begin{align}
& \mathcal{E}^L(t) +\int_0^t \mathcal{D}^L(s)\mm{d}{s} \leqslant c_2^2
e^{c_2 \sqrt{I_H^0}}I_L^0,\label{20210asfdsa1242212}\\
&{e}^{c_1 t}\mathcal{E}^H(t) + \int_0^t e^{c_1 \tau}\mathcal{D}^H(s)\mm{d}{s}
\leqslant c_2^2 I_H^0 .\label{estemsadfaalas}
\end{align}
\subsection{Proof of Theorem \ref{thm2} }\label{subsec:08}
Now we state the local well-posedness result for the initial-boundary value problem \eqref{01dsaf16asdfasf}.
\begin{pro}\label{202102182115}
Let $b >0$ be a constant and $\iota >0$ be the  constant in Lemma \ref{pro:1221}.
We assume that   $(  u^0,\eta^0)\in H^{1,3}_{\mm{s}}\times (H^{2,3}_{ \mm{s}}\cap H^3_*) $,
${  \|(u^0,\partial_1u^0,\partial_1^2\eta^0)\|_{3}\leqslant b }$  and $\mm{div}_{\mathcal{A}^0}u^0=0$, where $\mathcal{A}^0$ is defined by $\zeta^0$ and $ \zeta^0 = \eta^0+y$.
Then there exist a constant $\delta_2 \in (0,\iota/2]$  such that  if $\eta^0$   satisfies
\begin{align}
\| \eta^0\|_{\underline{1},3}\leqslant \delta_2,
\nonumber
\end{align}
the initial-boundary value problem  \eqref{01dsaf16asdfasf} admits a unique local classical solution
$(u,\eta,q)\in \mathfrak{U}_{T }^{1,3}\times   \mathfrak{C}^0 (\overline{I_{T }},{H}^{2,3}_{\mm{s}} )\times   C(\overline{I_T},\underline{H}^3)$  for some local existence time $T>0$.
Moreover, $(u,\eta,q) $ satisfies
\begin{align}& \sup\nolimits_{t\in \overline{I_T}}
\| \eta\|_{ 3}\leqslant 2\delta_2
,\nonumber \\
&\label{202012212151}
  \sup\nolimits_{t\in \overline{I_T}}\| (u,\eta, \partial_1 \eta )\|_3
  \leqslant c_3 \sqrt{I^0_L}  \mbox{ and }
  \sup\nolimits_{t\in \overline{I_T}}\| ( u , \partial_1 \eta)\|_{\underline{1},3}
  \leqslant c_3 \sqrt{ I^0_H}
\end{align}
for some positive constant $c_3\geqslant 1$.
It should be noted that $\delta_2$ and $T$ depend on $a$, $m $ and $\Omega$; moreover $T$ further depends on $b$.
\end{pro}
\begin{rem}
 {Here the uniqueness means that if there is another solution
$(\tilde{u}, \tilde{\eta},\tilde{q})\in\mathfrak{U}_{T }^{1,3}\times   \mathfrak{C}^0 (\overline{I_{T }},{H}^{2,3}_{\mm{s}} )\times   C(\overline{I_T},\underline{H}^3)$
 satisfying $0<\inf_{(y,t)\in \mathbb{R}^2\times \overline{I_T}} \det(\nabla \tilde{\eta}+I)$ and $\nabla_{(\nabla \tilde{\eta}+I)^{-\mm{T}}}\tilde{q}\in L^\infty_TH^3$, then
 $(\tilde{\eta},\tilde{u},\tilde{q})=(\eta,u,q)$ by virtue of the smallness condition
 ``$\sup_{t\in \overline{I_T}}\|\eta\|_{3}\leqslant 2\delta_2$''. In addition,
  we have, by the fact ``$\sup\nolimits_{t\in \overline{I_T}}\| \eta\|_3\leqslant  \iota$" and  Lemma \ref{pro:1221},
$$\inf\nolimits_{(y,t)\in \mathbb{R}^2\times  {I_T}} \det(\nabla \eta+I)\geqslant 1/4.$$}
\end{rem}
\begin{pf} The proof of Proposition \ref{202102182115} will be provided in Section \ref{202102241211}.
\hfill $\Box$
\end{pf}

Thanks to the \emph{priori} estimates \eqref{20210asfdsa1242212}, \eqref{estemsadfaalas}
and Proposition  \ref{202102182115}, we can easily establish
Theorem \ref{thm2}. Next we briefly give the proof.

Let  $(u^0,\eta^0)$ satisfies the assumptions in Theorem \ref{thm2},
\begin{align}
 e^{c_2 c_3\sqrt{I_H^0}}I_L^0  \leqslant\delta^2 \mbox{ and }
\delta =\min\left\{ \delta_1, \delta_2\right\}/ {{c}_2}c_3^2.\nonumber
\end{align}  By Proposition  \ref{202102182115}, there exists a unique local solution $( u,\eta,q)$ to the initial-boundary value  problem \eqref{01dsaf16asdfasf} with a maximal existence time $T^{\max}$, which satisfies that
\begin{itemize}
  \item for any $\tau\in  I_{T^{\max}}$,
the solution $(u,\eta,q)$ belongs to $\mathfrak{U}_{\tau}^{1,3}  \times{\mathfrak{H}}^{2,3}_{1,*,\tau}\times  C^0(\overline{I_\tau},\underline{H}^3)$, $$
 \sup\nolimits_{t\in {I_{\tau}}}\| \eta\|_{ 3}\leqslant 2\delta_2;$$
\item $\limsup_{t\to T^{\max} }\|\eta ( t)\|_{\underline{1},3} > \delta_2$
or $\limsup_{t\to T^{\max} }\|(u,\partial_1u,\partial_1\eta)( t)\|_{ 3}=\infty$, if $T^{\max}<\infty$.
\end{itemize}

Let
\begin{equation}
\nonumber
T^{*}:=\sup\left\{ \tau \in I_{T^{\max}}~\left|~
\|( u,\eta,\partial_1\eta)(t)\|_{3} \leqslant {{c}_2}c_3^2\delta  \;\mbox{ for any }t\leqslant \tau\right.\right\}.
\end{equation}
We easily see that the definition of $T^*$ makes sense and $T^*>0$.
Thus, to obtain the existence of global solutions,  it suffices to verify $T^*=\infty$. Next we shall prove this fact by contradiction.

We assume that $T^*<\infty$. Then, for any given $T^{**}\in I_{T^*}$,
\begin{equation}
\label{201911262sadf202}
\sup\nolimits_{t\in \overline{I_T^{**}}}\|(u,\eta,\partial_1\eta)(t) \|_3   \leqslant {{c}_2}c_3^2\delta \leqslant \delta_1.
\end{equation}
Thanks to the estimate  \eqref{201911262sadf202}, we can
follow the arguments of  \eqref{20210asfdsa1242212}--\eqref{estemsadfaalas} and use a regularity method and \eqref{202012212151} to derive that
\begin{align}
&\label{202104191546}
\|e^{c_1 t}\mathcal{E}^H(t)\|_{L^{\infty}(I_{T^{**}})} + \int_0^{{T^{**}}} e^{c_1 \tau}\mathcal{D}^H(s)\mm{d}s
\leqslant {c}_2^2
\lim_{t\to 0}\sup\nolimits_{\tau\in (0,t)}\|(u,\partial_1\eta)(\tau)\|_{\underline{1},3}^2
\leqslant c_2^2 c_3^2{ I^0_H},\\
&\label{202104191547}
\|\mathcal{E}^L(t)\|_{L^{\infty}(I_{T^{**}})} + \int_0^{T^{**}} \mathcal{D}^L(s)\mm{d}s
\leqslant {c}_2^2e^{c_2c_3\sqrt{I_{H}^0}}\lim_{t\to 0}\sup\nolimits_{\tau\in (0,t)}
 \|(u,\eta,\partial_1\eta)\|_3^2
\leqslant c_2^2c_3^2\delta^2.
\end{align}
By the weak continuity of $(\eta,u)$ and the fact
\begin{align}\sup\nolimits_{t\in \overline{I_{\tau}}} \|f\|_0
= \|f\|_{L^\infty(I_{\tau},L^2)}
\mbox{ for any }f\in C(\overline{I_{\tau}}, L^2_{\mm{weak}})\ \mbox{with}\ \tau>0, \nonumber
\end{align}
we further derive from \eqref{202104191546} and \eqref{202104191547} that
\begin{align}
& \sup\nolimits_{t\in \overline{I_T^{**}}}\|(u,\partial_1\eta)(t)\|_{\underline{1},3}
\leqslant c_2c_3 \sqrt{ I^0_H} \nonumber ,\\
&\label{202104191542}\sup\nolimits_{ t\in \overline{I_T^{**}}}\|(u,\eta,\partial_1\eta)(t)\|_{3}
\leqslant c_2c_3\delta.
\end{align}

We take $(u(T^{**},\eta(T^{**}) ))$ as a initial data. Noting that
$$\begin{aligned}
&\|(u,\partial_1u, \partial_1^2\eta)(T^{**})\|_{ 3}
\leqslant b := c_2c_3\sqrt{ I^0_H}
\mbox{ and }\| \eta(T^{**})\|_{\underline{1},3}\leqslant c_2c_3\delta\leqslant \delta_2,
\end{aligned}$$ then,
by Proposition \ref{202102182115} and \eqref{202104191542}, there exists a unique local  classical solution, denoted by $( u^*, \eta^*$, $q^*)$, to the initial-boundary value problem \eqref{01dsaf16asdfasf} with $(u(T^{**}),\eta(T^{**}) )$ in place of $(u^0,\eta^0 )$; moreover,
\begin{align}
&\sup\nolimits_{t\in [T^{**},T]}\| (u^* ,\eta^*, \partial_1 \eta^*  )\|_3 \leqslant {{c}_2}c_3^2\delta
\mbox{ and } \sup\nolimits_{t\in [T^{**},T]}\| \eta^*\|_{ 3}\leqslant 2\delta_2  , \nonumber
\end{align}
where the local existence time $T>0$ depends on $ b$, $a$, $m $ and $\Omega$.

In view of the existence result of $(u^*, \eta^*, q^*)$, the uniqueness conclusion in Proposition  \ref{202102182115}  and the fact that $T^{\max}$ denotes the maximal existence time, we immediately see that $T^{\max}>T^*+T/2$ and
$ \sup\nolimits_{t\in  {[0,T^*+T/2]}}\| (u,\eta, \partial_1 \eta)\|_3 \leqslant  {{c}_2}c_3^2\delta$.
This contradicts with the definition of $T^*$. Hence $T^*=\infty$ and thus $T^{\max}=\infty$.
This completes the proof of the existence of global solutions.
The uniqueness of global solution is obvious due to the uniqueness of local solutions in
Proposition \ref{202102182115} and the fact $\sup\nolimits_{t\geqslant0}\| \eta\|_{\underline{1},3}\leqslant 2\delta_2 $.
In addition,  it is obvious that $(\eta,u)$ indeed enjoys
\eqref{1.200} and \eqref{202101241112} by
 the derivation of \eqref{202104191546} and \eqref{202104191547}.

Finally, we derive \eqref{1.200xx}.
By \eqref{202101241112}, we easily see that
\begin{align}
\label{202101242112}
&\partial_1\eta(t)\to 0\mbox{ in }H^3\mbox{ as }t\to \infty
\end{align}
and
\begin{align}
\label{20safda2101242112}
\left\|\int_0^t u\mm{d}\tau\right\|_3\lesssim\int_0^t \|u  \|_3\mm{d}\tau\lesssim \|(u^0,\partial_1\eta^0)\|_{\underline{1},3}   \mbox{ for any }t>0\end{align}
Thanks to \eqref{20safda2101242112}, there exist a subsequence $\{t_n\}_{n=1}^\infty$ and some function $\eta^\infty \in H^3$ such that
$$\eta(t_n)-\eta^0= \int_0^{t_n} u\mm{d}\tau \to\eta^\infty -\eta^0\mbox{ weakly in }H^3.$$
Exploiting $\eqref{01dsaf16asdfasf}_1$ and \eqref{202101241112}, we have
$$
\begin{aligned}
\|\eta(t)-\eta^\infty\|_3\leqslant &\liminf_{t_n\to \infty} \left\|\int_t^{t_n } u\mm{d}\tau\right\|_3
\lesssim \sqrt{ I^0_H} \liminf_{t_n \to \infty}\int_t^{t_n} e^{-c_1 t/2}\mm{d}\tau
 \lesssim   \sqrt{ I^0_H}e^{-c_1 t/2}  ,
\end{aligned}$$
which, together with \eqref{202101242112}, yields that \eqref{1.200xx} holds and $\eta^\infty $ only depends on $y_2$ and $y_3$. This completes the proof of Theorem  \ref{thm2}.

\section{Local well-posedness}\label{202102241211}
 This section is devoted to the proof of local well-posedness of the initial-boundary value problem \eqref{01dsaf16asdfasf}.   Next we roughly sketch the proof frame to establish the  local well-posedness result.

Similarly to \cite{GXMWYJJMPA2019}, in which the authors investigated the well-posdeness problem of  incompressible inviscid MHD fluid equations with free boundary, we first  alternatively  use an iteration method to establish the existence result of the unique local solution $(u,\eta,{q})$ of the linear $\kappa$-approximate problem
\begin{equation} \label{01dsaf16asdfasf0000xx}
                              \begin{cases}
\eta_t-\kappa\partial_1^2\eta=u  ,\\
u_t+\nabla_{{\mathcal{B}}} {q}+a u=m   \partial_1^2\eta ,\\
\div_{{\mathcal{B}}} u=0  , \\
(u,\eta)|_{t=0}=(u^0,\eta^0) , \\
(u_3,\eta_3)|_{\partial\Omega}=0,\\
  \mathcal{B}_{3j}\partial_j{q}|_{\partial\Omega}=0,
\end{cases}
\end{equation}
where $\kappa>0$, $\mathcal{B}:=(\nabla \varsigma +I)^{-\mm{T}}$, and the given function $\varsigma$ belongs to
some proper function space,
see Proposition \ref{202102151544} in Section \ref{subsec:local01} for details. It should be noted that the Neumann boundary-value condition \eqref{01dsaf16asdfasf0000xx}$_6$ makes sure the solvability of $q$.

 Then we derive the $\kappa$-independent estimates of the solutions, denoted by $(u^\kappa,\eta^\kappa, {q}^\kappa)$ for $\kappa\in \mathbb{R}^+$, of the $\kappa$-approximate problem above,
and then  take the limit of $(u^\kappa,\eta^\kappa, {q}^\kappa)$ with respect to $\kappa\to 0$ in a common definition domain. In particular, the limit function, denoted by $(u,\eta,{q})$, is the unique local solution to the linearized problem
\begin{equation}\label{01dsaf16asdfsafasf0000}
                              \begin{cases}
\eta_t =u,\\
u_t+\nabla_{{\mathcal{B}}} {q}+a u=m   \partial_1^2\eta,\\
\div_{{\mathcal{B}}} u=0, \\
(u,\eta)|_{t=0}=(u^0,\eta^0), \\
(u_3,\eta_3)|_{\partial\Omega}=0,
\end{cases}
\end{equation}
see Proposition \ref{thm08} in Section \ref{subsubsce:03} for details. It should be noted that the solution $q$ of \eqref{01dsaf16asdfsafasf0000} automatically satisfies \eqref{01dsaf16asdfasf0000xx}$_6$ due to $u_3$ satisfying $\partial_t u_3 |_{\partial\Omega}=0$.

Finally, since the linearized problem \eqref{01dsaf16asdfsafasf0000} admits a unique local solution for any given $\varsigma$, thus we easily arrive at Proposition \ref{202102182115}  by a standard iteration method as in \cite{Coutand1,Coutand2}, in which the authors investigated the well-posdeness problem of  incompressible Euler equations with free boundary.

Now we turn to introducing some new notations appearing in this section.
$P(x_1,\ldots, x_n )$
and $\dot{P}(x_1,\ldots, x_n )$ represent the generic polynomials with respect to the parameters $x_1$, $\ldots$, $x_n$, where all the coefficients in $P$ and $\dot{P}$ are equal one, and $\dot{P}$ satisfies $\dot{P}(0,\ldots,0)=0$;
$a\lesssim_\kappa b$ means that $a\leqslant c^\kappa b$, where $c^\kappa$ denotes a generic positive constant, which may depend on $\kappa$, $a$, $m $ and $\Omega$. It should be noted that $P(x_1,\ldots, x_n )$, $\dot{P}(x_1,\ldots, x_n )$ and $c^\kappa$ may vary from line to line.

We define the difference quotient with respect variable $y_1$ as follows:
\begin{align}
\label{202220105041846}
D_1^{\tau}f (y_1,y_2,y_3)=(f(y_1+\tau,y_2 ,y_3)-f(y_1,y_2,y_3 ))/\tau\mbox{ for }
|\tau|\in  (0,1).
\end{align}
 The notation $\mathcal{X}$ represents $\mm{Id}$\mbox{ or }$D^\tau_1$, where $\mm{Id}$ denotes and identity mapping, and $D^\tau_1$ the difference quotient with respect to the variable $y_1$.

We always use the notations $\mathcal{B}$, resp. $\mathcal{J}$ to represent $(\nabla \varsigma +I)^{-\mm{T}}$, resp. $\det(\nabla \varsigma+I)$, where $\varsigma$ at least satisfies
\begin{align}
\varsigma \in {C}^0(\overline{I_T},{H}^{3})\mbox{ and }\inf\nolimits_{(y,t)\in {\Omega_T}}\det(\nabla \varsigma +I)\geqslant 1/4 \mbox{ for some }T>0 .  \label{202safas2103250842}
\end{align}
In addition, we define that  \begin{align}
\mathbb{A}^{3,1/4}_{T,\iota }:=\{\psi\in  {C}^0(\overline{I_T},{H}^{3}_{\mm{s}})~|~  \|\psi\|_3\leqslant \iota\}  \label{2022104161633}
  \end{align}
and
\begin{align}
\mathbb{S}_T:=\{&(w,\xi,\beta) \in C^0(\overline{I_{T}} ,{H}^3_{ \mathrm{s}})
\times C^0( \overline{I_{T}},{H}^{3}_{\mm{s}})\times (C ^0(\overline{I_{T}},\underline{H}^3 )  )~|~   \partial_1^2\xi  \in C^0(\overline{I_T},H^2_{\mm{s}}),\nonumber \\
& \partial_1^2\xi \in  L^2_{T}H^{3} ,\
\nabla^3  \partial_1\xi
\in C^0_{ B, \mm{weak}}( \overline{I_T},L^2), \  \nabla_{{\mathcal{B}}} \beta\in {L_{T}^2H^{ 3}}, \;
 {{\mathcal{B}}}_{3j}\partial_j \beta
|_{\partial\Omega}=0  \}  \label{20210413asfda21251},
\end{align}
 where $T>0$ and $\iota$ is the positive constant provided in Lemma \ref{pro:1221}. \emph{It should be noted that the function $\varsigma$, which belongs to $\mathbb{A}^{3,1/4}_{T,\iota }$, automatically satisfies \eqref{202safas2103250842} by Lemma \ref{pro:1221}. }

Finally, some preliminary estimates for $\mathcal{B}$ and $\mathcal{J}$ are collected as follows.
\begin{lem}\label{lem:0817}
Let $\varsigma $ satisfy \eqref{202safas2103250842} and $i=0$, $1$. Then
\begin{align}
&\label{202008121505}
\| \mathcal{B} \|_{C^0(\bar{\Omega})}
+\| (\mathcal{B},\mathcal{J}^{-1}) \|_2
\lesssim  P(\|{\varsigma}\|_3),\\
&\label{202008121535nnn}
\|{\mathcal{B}}-I\|_{ {i},2}
\lesssim  \|{\varsigma}\|_{{i},3}{P}(( \|{\varsigma}\|_{3} ) \mbox{ for any }t\in \overline{I_T}.
\end{align} If additionally  $\varsigma$ further satisfies
$(\partial_1\varsigma,\varsigma_t)\in  L^\infty_TH^{i,3}$ , then, for a.e. $t\in {I_T}$,
\begin{align}
&\label{202008121535}
\| \partial_t(\mathcal{B},\mathcal{J})\|_{i, 2}
\lesssim
          \dot{P}( \|( \varsigma ,\varsigma _t)\|_{\underline{i}, 3})  ,
\\
&\label{202008121535n}
\|  \mathcal{B}\|_{ 1+i, 3}
\lesssim
      \dot{P}( \| {\varsigma}\|_{\underline{1+i},3}).
\end{align}
\end{lem}
\begin{pf} Since the derivation of the estimates
\eqref{202008121505}--\eqref{202008121535n} are very elementary, we omit it. \hfill$\Box$
\end{pf}

\subsection{Solvability of the linear $\kappa$-approximate problem \eqref{01dsaf16asdfasf0000xx}}\label{subsec:local01}

In this section, we construct the unique local solution of the  linear $\kappa$-approximate problem \eqref{01dsaf16asdfasf0000xx}. To this purpose, we shall rewrite  the linear $\kappa$-approximate problem \eqref{01dsaf16asdfasf0000xx} as the following equivalent problem (in the sense of classical solutions):
\begin{equation} \label{01dsaf16asdfasf0safd000}
                              \begin{cases}
\eta_t-\kappa\partial_1^2\eta=u ,\\
u_t+\nabla_{{\mathcal{B}}} {q}+a u=m  \partial_1^2\eta  ,\\
-\Delta_{{\mathcal{B}}}{q}=K^1  , \\
(u,\eta)|_{t=0}=(u^0,\eta^0)  , \\
(u_3,\eta_3)|_{\partial\Omega}=0 ,\ \mathcal{B}_{3j}\partial_j{q}|_{\partial\Omega}=0,
\end{cases}
\end{equation}
where $\Delta_{{\mathcal{B}}}:=\mm{div}_{{\mathcal{B}}}\nabla_{{\mathcal{B}}}$, $\varsigma|_{t=0}=\eta^0$, $(u^0,\eta^0)$ satisfies
$\mm{div}_{(\nabla \eta^0+I)^{-\mm{T}}}u^0=0$, $\inf\nolimits_{y \in \overline{\Omega}} \det{(\nabla \eta^0+I)^{-\mm{T}}} >0$, and
\begin{align}
& K^1:=a\mm{div}_{{\mathcal{B}}}{u} -\mm{div}_{\mathcal{B}_t}u-\mathcal{J}_t\mm{div}_{{\mathcal{B}}}u/\mathcal{J}
  -m   \mm{div}_{{\mathcal{B}}} \partial_1^2{\eta}.  \label{202102141605}
\end{align}
Then the solvability of the linear $\kappa$-approximate problem  \eqref{01dsaf16asdfasf0000xx}
reduces to the solvability of the initial-boundary value problem  \eqref{01dsaf16asdfasf0safd000}.

We want to establish the local well-posdeness result for \eqref{01dsaf16asdfasf0safd000}
by an iteration method.
To begin with, we shall investigate the solvability of the linear problem
\begin{equation} \label{01dsa000}
                              \begin{cases}
\eta_t-\kappa\partial_1^2\eta=w , \\
u_t+a u= K^2, \\
-\Delta_{{\mathcal{B}}}{q}=K^1 ,  \\
(u,\eta)|_{t=0}=(u^0,\eta^0) ,   \\
(u_3,\eta_3)|_{\partial\Omega}=0 ,\ \mathcal{B}_{3j}\partial_j{q}|_{\partial\Omega}=0,
\end{cases}
\end{equation}
where $(\xi,w,\theta)\in \mathbb{S}_T$ is given, and $K^2:=m    \partial_1^2\xi-\nabla_{{\mathcal{B}}}\theta  $.

It is easy to see that the above linear problem is  equivalent to the following three sub-problems: the  initial-boundary value problem of partly dissipative equation for $\eta$
\begin{equation}\label{parabolic}
                              \begin{cases}
\eta_t-\kappa\partial_1^2\eta=w ,\\
\eta|_{t=0}=\eta^0  ,\\
\eta_3|_{\partial\Omega} =0  ,
\end{cases}
\end{equation}
the initial-value problem of ODE  for $u$
\begin{equation}\label{transport}
                              \begin{cases}
 u_t+a u=K^2 ,\\
u|_{t=0}=u^0  , \\
u_3|_{\partial\Omega}=0
\end{cases}
\end{equation}
and the Neumann boundary-value problem of elliptic equation for $q$
\begin{equation}\label{elliptic}
                              \begin{cases}
-\Delta_{{\mathcal{B}}}{q}=K^1  &\mbox{in } \Omega,\\
\mathcal{B}_{3j}\partial_j{q} =0 &\mbox{on } \partial\Omega.
\end{cases}
\end{equation}
Thus the solvability of the linear problem \eqref{01dsa000} reduces to the solvability of the  three sub-problems above. Next we establish the { global} well-posedness results for the above  three sub-problems  in sequence.
 \begin{pro}\label{pro:parabolic}
Let $\kappa>0$, $T>0$, $w\in L^2_{T}H^{3}_{\mm{s}}$
and $\eta^0\in H^{1,3}_{\mm{s}}$.
Then the initial-boundary value problem \eqref{parabolic} admits a unique solution $\eta\in {C}^0( \overline{I_T},{H}^{3}_{\mm{s}} )$, which satisfies $ \partial_1^2 \eta \in C^0(\overline{I_T},H^2_{\mm{s}})$,
$\nabla^{3}\partial_1 \eta \in C^0_{B,\mm{weak}}(\overline{I_T},{L}^2)$,
and enjoys the estimates
\begin{align}
&
\sup\nolimits_{t\in \overline{I_T}}(\| \eta \|_{ 3}+\sqrt{\kappa}\| \eta \|_{ 1,3})
+ \|(\kappa\partial_1^{2 }\eta,\eta_t)\|_{L^2_TH^{3}}\nonumber \\
&\lesssim \|\eta^0 \|_{ 3}+\sqrt{\kappa}\|\eta^0 \|_{ 1,3}+(1+\sqrt{T})\|  w\|_{L^2_TH^{3}}, \label{20200805100sadfsa5}\\
&\label{20200805DSAFS100sadfsa5}\ {\kappa}\sup\nolimits_{t\in\overline{I_s}} \|   \eta^\varepsilon  \|_{1,3}^2
 \lesssim  \kappa\| \eta^0 \|_{1 ,3}^2+  \| w\|_{{L^2_sH^{3}}}^2\mbox{ for any }s\in (0,T].
\end{align}
\end{pro}
\begin{pf}
Let $\varepsilon\in (0,1)$,
$\chi$ be a one-dimensional standard mollifier (see \cite[pp. 38]{NASII04} for the definition), and $\chi^{\varepsilon}(s):=\chi(s/\varepsilon)/\varepsilon$.
Let $\tilde{w}=w$ in $\Omega_T $ and $\tilde{w}=0$ outside $\Omega\times (\mathbb{R}\backslash I_T)$. We define the mollification of $\tilde{w}$ with respect to $t$:
\begin{align}
\nonumber
S^{t}_\varepsilon (\tilde{w}): =\chi^{\varepsilon} * \tilde{w}.
\end{align}
Then
$S^{t}_\varepsilon  (\tilde{w}) \in C^\infty(\mathbb{R},H^{1,3}_{\mm{s}})$. We can check that
\begin{align}
&\|  S^{t}_\varepsilon (\tilde{w} ) \|_{L^2_TH^{3}}  \lesssim \| w\|_{L^2_TH^{3}}
\label{202210208asdfsa2248},\\
&  S^{t}_\varepsilon  (\tilde{w}) \to  w
\mbox{ strongly in }{L^2_TH^{3}},
\nonumber
\end{align}

Now we consider the  $\tau$-approximate problem for \eqref{parabolic}:
\begin{equation}\label{parabsalic}
                              \begin{cases}
 \eta^\tau_t= {L}(\eta^\tau)+ S^{t}_\varepsilon  (\tilde{w} )
 ,\\
\eta^\tau|_{t=0}=\eta^0 ,\\
\eta ^\tau_3|_{\partial\Omega}=0 ,
\end{cases}
\end{equation}
where $L:$ $H^{1,3}_{\mm{s}}\to H^{1,3}_{\mm{s}}$
by the rule $L(f)= \kappa D^{-\tau}_1 (D^\tau_1f)$ for $f\in  H^{1,3}_{\mm{s}}$, see  \eqref{202220105041846} for the definition of difference quotient $D^\tau_1$.

It is easy to see  that  $L(\varpi)\in  C^0(\overline{I_T},H^{1,3}_{\mm{s}})$ for $\varpi\in  C^0(\overline{I_T},H^{1,3}_{\mm{s}})$ and $L\in \mathcal{L}(H^{1,3}_{\mm{s}})$,
where $\mathcal{L}(H^{1,3}_{\mm{s}})$ is a set of all linear bounded operators of $H^{1,3}_{\mm{s}}$.
In particular, $ L\in  C^0(\overline{I_T},\mathcal{L}(H^{1,3}_{\mm{s}}))$.
By existence theory of the initial-value problem of a abstract ODE equation
(see \cite[Proposition 2.17]{NASII04}),
there exists a unique solution
$\eta^\tau\in C^0(\overline{I_T},H^{1,3}_{\mm{s}})\cap C^1( \overline{I_T},H^{1,3}_{\mm{s}})$ to \eqref{parabsalic}.
Obviously
$\eta^\tau$,  $\eta^\tau_t$ automatically belong to $L^2_TH^{1,3}_{\mm{s}}$
by the second conclusion in Lemma \ref{20021032019018} and the separability of $H^{1,3}_{\mm{s}}$.

Let $\alpha$ satisfy $0\leqslant |\alpha|\leqslant  3$.
Applying $\partial^\alpha $
to \eqref{parabsalic}$_1$, and then multiplying the resulting identity by $ \partial^\alpha  \eta^\tau$  in $L^2 $, we have
\begin{align}
 \frac{\mm{d}}{\mm{d}t} \int |\partial^{\alpha} \eta^\tau|^2\mm{d}y
+ \kappa\int |D_1^{\tau}\partial^{\alpha} \eta^\tau |^2\mm{d}y
  =\int S^{t}_\varepsilon ( \partial^{\alpha} \tilde{w}) \cdot
\partial^{\alpha} \eta^\tau \mm{d}y . \label{202008051043}
\end{align}
Making use of  \eqref{202210208asdfsa2248},
we easily deduce from \eqref{202008051043} that,
\begin{align}
\| \eta^\tau\|_{C^0(\overline{I_T},H^{3})}
\lesssim\| \eta^0 \|_{ 3}+  \sqrt{T}\| w\|_{L^2_TH^{3}}.  \label{202102082127}
\end{align}

In addition, we easily deduce from  \eqref{parabsalic}$_1$ that
\begin{align}\label{2020safda08051043}
& \frac{\kappa}{2}\frac{\mm{d}}{\mm{d}t}\| D_1^{\tau}\partial^{\alpha} \eta^\tau\|^2_0
+\int| \partial^{\alpha} \eta^\tau_t |^2\mm{d}y
=\int (\partial^{\alpha} S^{t}_\varepsilon (\tilde{w}))
\cdot  \partial^{\alpha} \eta^\tau_t \mm{d}y.
\end{align}
Noting that
\begin{align}
\label{2021032101049}
\| D_1^{\tau}\partial^{\alpha} \eta^\tau|_{t=0}\| _0
\lesssim \| \partial^{\alpha}\eta^0\| _{1 ,0},
\end{align}
thus, making use of  Young's inequality, \eqref{202210208asdfsa2248}, \eqref{parabsalic}$_1$ and \eqref{2021032101049}, we deduce from \eqref{2020safda08051043} that, for any $s\in (0,T] $
\begin{align}
 & {\kappa}\sup\nolimits_{t\in\overline{I_s}} \| D_1^\tau \eta^\tau \|_{3}^2
+ \| ({L}( \eta^\tau),  \eta^\tau_t)\|_{L^2_sH^{3}}^2
\lesssim  \kappa\| \eta^0 \|_{1 ,3}^2+  \| w\|_{{L^2_sH^{3}}}^2.  \label{202008051056}
\end{align}

Thanks to the regularity of $\eta^\tau$, we can we easily derive \eqref{parabsalic}$_1$ that, for any $\varphi\in H^1$ and for any $r$, $s\in I_T$,
\begin{align}
& \int D_1^\tau\partial^{\alpha}
(\eta^\tau(y,t) -  \eta^\tau(y,r))\cdot\varphi\mm{d}y \nonumber \\
&=
- \int_r^t\int ( S^{t}_\varepsilon ( \partial^{\alpha}  \tilde{w})   +
L(\partial^{\alpha}  \eta^\tau)) \cdot D_1^{-\tau}\varphi  \mm{d}y \mm{d}s.
\label{202210041111912}
\end{align}
Exploiting  \eqref{202210208asdfsa2248},
the uniform estimate of $L( \eta^\tau)$ in \eqref{202008051056} and the fact
$$\|D_1^{-\tau}\varphi\|_0\lesssim \| \varphi\|_{1,0} ,$$
we easily derive from the identity \eqref{202210041111912} that
\begin{align}
&\label{20202108021asdf27}
 D_1^{\tau}\nabla^3  \eta^\tau\mbox{ is uniformly continuous in }H^{-1} .
\end{align}
Here and in what follows, $H^{-1}$ denotes the dual space of $H^1_0:= \{w\in {H}^1~|~w|_{\partial\Omega} =0\} $.

Making use of \eqref{202102082127}, \eqref{202008051056}, \eqref{20202108021asdf27} and \eqref{202104141653asfda}--\eqref{202104141653} there exists a sequence of $\{\eta^\tau\}_{|\tau|\in (0,1)}$ (still denoted by $\eta^\tau$ for simplicity) such that, as $\tau\to 0$,
\begin{align}
&( \eta_t^\tau , D_1^{-\tau}(D_1^\tau \eta^\tau) )\rightharpoonup  (\eta_t^\varepsilon ,\partial_1^{2 } \eta^\varepsilon )\mbox{ weakly in }L^2_TH^{3}_{\mm{s}}, \
 \eta^\tau \to  \eta^\varepsilon  \mbox{ strongly in }
C^0(\overline{I_T},H^{2}), \nonumber  \\
& (\eta^\tau,D_1^\tau \eta^\tau ) \rightharpoonup (\eta^\varepsilon,\partial_1 \eta^\varepsilon) \mbox{ weakly-* in }L^\infty_TH^{3}_{\mm{s}} ,\
  D_1^\tau\nabla^3 \eta^\tau \to\nabla^3\partial_1 \eta^\varepsilon
 \mbox{ in } C^0(\overline{I_T},L^2_{\mm{weak}}) . \nonumber
\end{align}
Moreover, the limit function $\eta^\varepsilon$ is just the unique solution to the problem \begin{equation}
\nonumber
                              \begin{cases}
 \eta^\varepsilon_t=  \kappa\partial_1^2  \eta^\varepsilon + S^{t}_\varepsilon  (\tilde{w} )
  ,\\
\eta^\varepsilon|_{t=0}=\eta^0  ,\\
\eta ^\varepsilon_3|_{\partial\Omega}=0
\end{cases}
\end{equation} and satisfies the following estimates:
\begin{align}
&\sup\nolimits_{t\in \overline{I_T}} \| \eta^\varepsilon  \|_{ 3}
 \lesssim \|\eta^0 \|_{ 3} +  \sqrt{T} \|  w\|_{L^2_TH^{3}} , \label{20200805100sasfsdafadfsa5} \\
& {\kappa}\sup\nolimits_{t\in\overline{I_s}} \|  \eta^\varepsilon  \|_{1,3}^2
+ \|(\kappa\partial_1^2\eta^\varepsilon,\eta^\varepsilon_t )\|_{L^2_sH^{3}}^2
\lesssim  \kappa\| \eta^0 \|_{1 ,3}^2+  \| w\|_{{L^2_sH^{3}}}^2\mbox{ for any }s\in (0,T].\label{202fsa5}
\end{align}
Thus we further have $ \partial_1^2 \eta^\varepsilon \in C^0(\overline{I_T},H^2_{\mm{s}})$ and  $\eta^\varepsilon \in C^0(\overline{I_T}, H^{ 3}_{\mm{s}})$.

Noting that $\eta^{\varepsilon}$ enjoys the estimates \eqref{20200805100sasfsdafadfsa5} and \eqref{202fsa5}, thus we  get a limit function $\eta$, which satisfies the desired conclusion in Proposition \ref{pro:parabolic}, by taking limit of some sequence of $\{\eta^\varepsilon\}_{\varepsilon>0}$. We omit such limit process, since it is very similar to the argument of obtaining $\eta^\varepsilon$.
\hfill $\Box$
\end{pf}

\begin{pro}\label{pro:transport}
Let $a\in \mathbb{R}$, $K^2\in L^1_{T}H^{3}_{\mm{s}}$ and $u^0\in H^{3}_{\mm{s}} $, then the initial-boundary value problem \eqref{transport} admits a unique solution $u\in C^0(\overline{I_T},H^{3}_{\mm{s}})$, which satisfies
$u_t\in L^1_{T}H^{3}_{\mm{s}}$ and
\begin{align}
&\label{202104101021}
\|u(t)\|_{ 3}\lesssim\|u^0\|_{ 3}
+\int_0^t\|K^2(s)\|_{ 3}\mm{d}s.
\end{align}
\end{pro}
\begin{pf}
Let
\begin{align}\label{202008051415}
u=u^0e^{-at}+\int_0^tK^2(s)e^{-a(t-\tau)}\mm{d}s.
\end{align}
Since $K^2\in L^1_{T}H^{ 3}_{\mm{s}}$ and $u^0\in H^{  3}_{\mm{s}}$, it is easy to check that $u$ given by \eqref{202008051415} is the unique solution of \eqref{transport};
moreover, $u$ belongs to $C^0(\overline{I_T},H^{ 3}_{\mm{s}})$ and satisfies \eqref{202104101021}.
\hfill $\Box$
\end{pf}

\begin{pro}\label{pro:elliptic}
Let $T>0$,  $\iota$ be the positive constant provided in Lemma \ref{pro:1221}
and $\varsigma\in \mathbb{A}^{3,1/4}_{T,\iota }$  defined by \eqref{2022104161633}.
If $K^1\in C^0(\overline{I_T},H^1)$ and satisfies $\int  K^1\mathcal{J}\mm{d}y=0$ for $t\in \overline{I_T}$, then there exist a unique solution ${q}\in C^0(\overline{I_T},\underline{H}^3)  $, which satisfies the boundary-value problem \eqref{elliptic} for each $t\in I_T$ and enjoys the estimate
\begin{align}
\label{202008051safa442}
&\| q\|_{C^0(\overline{I_T}, {H}^3) }\lesssim P( \|{\varsigma}\|_{L^\infty_TH^{3}}) \|K^1\|_{L^\infty_TH^{1}} .
\end{align}

If we additionally assume that $K^1\in   L^r_TH^2$ for $r=2$ or $\infty$, then
\begin{align}\label{202008051442}
&\|\nabla_{{\mathcal{B}}}q\|_{L^r_TH^{3}}\lesssim P( \|{\varsigma}\|_{L^\infty_TH^{ 3}}) \|K^1\|_{L^r_TH^2} .
\end{align}
\end{pro}
\begin{pf}
 Let $\varphi=\varsigma+y $, then $\varphi$ is a $C^1$-diffeomorphism mapping by Lemma \ref{pro:1221}.
We denote $\varphi^{-1}$ the inverse mapping of $\varphi$ with respect to variable $y$, and then define $\tilde{K}^1:=K^1 (\varphi^{-1},t)$. Obviously $\tilde{K}^1\in  C^0(\overline{I_T}, H^{1})$ by \eqref{202104031901}.

Now we consider the following Neumann boundary-value problem of elliptic equation
\begin{align}\label{elliptic1}
\begin{cases}
-\Delta p=\tilde{K}^1&\mbox{in } \Omega,\\
\partial_3 p =0 &\mbox{on } \partial\Omega.
\end{cases}
\end{align}
It is easy to see that
\begin{align}
\int \tilde{K}^1\mm{d}x=\int K^1\mathcal{J}\mm{d}y=0 .
\nonumber
\end{align}

By Lemma \ref{lem:08181945}, there exists
a unique solution $p\in C^0(\overline{I_T},\underline{H}^3) $ of the boundary-value problem \eqref{elliptic1} such that
\begin{align}\label{2020080516sasafddfa42}
\|p \|_{ C^0(\overline{I_T}, {H}^3)}\lesssim \|\tilde{K}^1\|_{L^\infty_TH^{1}}.
\end{align}

Let $\tilde{q}=p(\varphi)$, then $\tilde{q}\in C^0(\overline{I_T}, {H}^3) $
by \eqref{2020103250855}. We further define that $ {q}:=\tilde{q}- (\tilde{q})_{\Omega} $, then $ {q} \in C^0(\overline{I_T}, \underline{H}^3)$ and satisfies  \eqref{202008051safa442} by \eqref{2020080516sasafddfa42} and \eqref{2022104101908}. Since $p$ satisfies \eqref{elliptic1}, we see that $q$ is the solution of the boundary-value problem \eqref{elliptic}. In addition, the uniqueness of the solution of the problem \eqref{elliptic} is obvious in the class $C^0(\overline{I_T}, \underline{H}^3)$.

If we additionally assume that $K^1\in   L^r_TH^{2}$, then  $\tilde{K}^1\in  L^r_TH^{2}$, where $r=2$ and $\infty$. By the second assertion in Lemmas \ref{20021032019018}--\ref{lem:08181945}, \eqref{202104032134} and the separability of $\underline{H}^3$,   we have $p\in L^r_TH^{3}$ and \begin{align}\label{2020080516sadfa42}
\|p \|_{L^r_TH^{4}}\lesssim P(\|{\varsigma}\|_{L^\infty_TH^{3}}) \|\tilde{K}^1\|_{L^r_TH^{2}}.
\end{align}
By \eqref{2020080516sadfa42}, \eqref{2021sfa04031901} and \eqref{2022104101908},
\begin{align}
&\|\nabla_{\mathcal{B}}q\|_{L^r_TH^{3}}=\|\nabla p(x,t)|_{x=\varphi(y,t)}\|_{L^r_TH^{3}}\nonumber \\
&\lesssim P(\|\varsigma\|_{L^\infty_TH^{3}}) \|\nabla p\|_{L^r_TH^{3}}\lesssim P(\|{\varsigma}\|_{L^\infty_TH^{3}}) \|{K}^1\|_{L^r_TH^{2}}, \nonumber
\end{align}
which yields \eqref{202008051442}.  \hfill $\Box$
\end{pf}

With Propositions \ref{pro:parabolic}--\ref{pro:elliptic} in hand, next we will use an
iterative method  to establish the existence result of a unique local solution to the linear $\kappa$-approximate problem \eqref{01dsaf16asdfasf0000xx}.

\begin{pro}
\label{202102151544}
Let $\mathbb{A}^{3,1/4}_{\alpha ,\iota }$, resp.  $\mathbb{S}_{ \alpha  }$ be defined by \eqref{2022104161633} resp. \eqref{20210413asfda21251} with some positive constant $\alpha  $ in place of $T$.
 We assume that  $\kappa>0$,  $a\geqslant 0$,  $(u^0,\eta^0)\in  {H}^{3}_{\mm{s}}\times{H}^{1,3}_{\mm{s}} $  and $\varsigma\in\mathbb{A}^{3,1/4}_{\alpha ,\iota }$  satisfies  $ \varsigma_t \in C^0(\overline{I_{\alpha }},H^{2})\cap {L^{\infty}_{\alpha }H^{3}}$, then  the $\kappa$-approximate problem \eqref{01dsaf16asdfasf0000xx} defined on $\Omega_\alpha$  admits a unique solution, denoted by $(u^\kappa,\eta^\kappa,q^\kappa)$, which belongs to $\mathbb{S}_{ \alpha  }$.
\end{pro}
\begin{pf}
Let $T  \leqslant \alpha $. In view of Propositions \ref{pro:parabolic}--\ref{pro:elliptic}, we easily see that there exists  a solution sequence $\{ (u^n,\eta^n,q^n)\}_{n=1}^\infty$ defined on $\overline{I_T}$. Moreover,
\begin{enumerate}
\item[(1)]  $(u^1,\eta^1, q^1)=0$;
\item[(2)]
$(u^n,\eta^n,Q^n)\in \mathbb{S}_{T  }$ for $n\geqslant 1$;
\item[(3)]
$(u^n,\eta^n,q^n)$ for $n\geqslant2$ enjoys the following relations:
\begin{equation}
\nonumber
                              \begin{cases}
{\eta}_t^{n}-\kappa\partial_1^2{\eta}^{n}={u}^{n-1},\\
{u}^{n}_t+a {u}^{n}= m   \partial_1^2{\eta}^{n-1}
-\nabla_{{\mathcal{B}}}{q}^{n-1}=:K^{2,n},\\
-\Delta_{{\mathcal{B}}}{q}^{n}
=a\mm{div}_{{\mathcal{B}}}{u}^{n}-\mm{div}_{\mathcal{B}_t}{u}^{n}
-\mathcal{J}_t\mm{div}_{{\mathcal{B}}}{u}^{n}/\mathcal{J}
-m  \mm{div}_{{\mathcal{B}}} \partial_1^2{\eta}^{n} =:K^{1,n}  , \\
({u}^{n},{\eta}^{n})|_{t=0}=(u^0,\eta^0) ,\\
( {u}^{n}_3,{\eta}^{n}_3)|_{\partial\Omega}=0 ,\
{\mathcal{B}}_{3j}\partial_j{q}^{n} |_{\partial\Omega} =0,
\end{cases}
\end{equation}
where $K^{1,n} \in  C^0(\overline{I_T  },H^1)\cap L^2_T  {H}^{2} $,
$K^{2,n} \in  L^2_{T  }H^{ 3}_{\mm{s}}$ and
\begin{align}
 \int  K^{1,n}\mathcal{J}\mm{d}y=0\mbox{ for each }t\in \overline{I_T} . \nonumber
\end{align}
It should be noted that we have used the relation $\partial_{j}(\mathcal{J} \mathcal{B}_{ij})=0$ for $i=1$, $2$, and the regularity $ \varsigma_t \in C^0(\overline{I_{\alpha }},H^1)$ in the derivation of the above identity.
\end{enumerate}

 We  further define
$(\bar{u}^{n+1},\bar{\eta}^{n+1}, \bar{q}^{n+1}):=(u^{n+1}-u^{n},\eta^{n+1}-\eta^{n},q^{n+1}-q^{n})$
and $(\bar{u}^{n},\bar{\eta}^{n}, \bar{q}^{n}):=(u^{n}-u^{n-1},\eta^{n}-\eta^{n-1},q^{n}-q^{n-1})$ for $n\geqslant2$.
Then we have
\begin{equation}\label{01dsaf16asdfasf00001}
                              \begin{cases}
\bar{\eta}_t^{n+1}-\kappa\partial_1^2\bar{\eta}^{n+1}=\bar{u}^{n},\\
\bar{u}^{n+1}_t+a \bar{u}^{n+1}= m   \partial_1^2\bar{\eta}^{n}
-\nabla_{{\mathcal{B}}}\bar{q}^{n},\\
-\Delta_{{\mathcal{B}}}\bar{q}^{n+1}
=a\mm{div}_{{\mathcal{B}}}\bar{u}^{n+1}-\mm{div}_{\mathcal{B}_t}\bar{u}^{n+1}
-\mathcal{J}_t\mm{div}_{{\mathcal{B}}}\bar{u}^{n+1}/\mathcal{J}
-m  \mm{div}_{{\mathcal{B}}}(\partial_1^2\bar{\eta}^{n+1})  , \\
(\bar{u}^{n+1},\bar{\eta}^{n+1})|_{t=0}=(0,0) ,\\
(\bar{u}^{n+1}_3,\bar{\eta}^{n+1}_3)|_{\partial\Omega}=0 ,\\
{\mathcal{B}}_{3j}\partial_{i}\bar{q}^{n+1} |_{\partial\Omega} =0
 .
\end{cases}
\end{equation}

Applying   Propositions \ref{pro:parabolic}--\ref{pro:elliptic} to \eqref{01dsaf16asdfasf00001}, and then using \eqref{202008121505} and \eqref{202008121535}, we get
\begin{align}
& \nonumber
\| \bar{\eta}^{n+1} \|_{C^0(\overline{I_T}, H^{3} )}
+\|\partial_1  \bar{\eta}^{n+1} \|_{L_{T}^{\infty}H^{3}}
+\| \partial_1^2\bar{\eta}^{n+1} \|_{L_{T}^2H^3}
\lesssim_{\kappa}T(1+\sqrt{T})\|\bar{u}^{n}\|_{L^\infty_TH^{3} },\\
&  \nonumber
\| \bar{u}^{n+1}  \|_{C^0(\overline{I_T}, H^{3} )}
\lesssim \sqrt{T} \|(\partial_1^2\bar{\eta}^{n},\nabla_{{\mathcal{B}}}\bar{q}^{n} )\|_{L_{T}^2H^{ 3}}, \\
&\|\bar{q}^{n+1}\|_{C^0(\overline{I_T}, H^{3} )}\lesssim_\kappa
P ( \|({\varsigma}, {\varsigma}_t)\|_{L^{\infty}_{T}H^3})
\|(\bar{u}^{n+1},\partial_1^2\bar{\eta}^{n+1})\|_{L_{T}^\infty H^2}, \nonumber  \\
&\|\nabla_{{\mathcal{B}}}\bar{q}^{n+1}\|_{L_{T}^2H^{ 3}}\lesssim_\kappa
P ( \|({\varsigma}, {\varsigma}_t)\|_{L^{\infty}_{T}H^{ 3}})
\left( \sqrt{T}\|\bar{u}^{n+1}\|_{L_{T}^{\infty}H^{ 3}}+\|\partial_1^2\bar{\eta}^{n+1}\|_{L_{T}^2H^{ 3}}\right). \nonumber
\end{align}
In addition, by \eqref{01dsaf16asdfasf00001}$_1$  and \eqref{01dsaf16asdfasf00001}$_2$,
$$\|(\bar{\eta}_t^{n+1},\bar{u}_t^{n+1})\|_{C^0(\overline{I_T}, H^{2} )}\lesssim \|( \bar{u}^{n+1},\bar{u}^{n},\partial_1^2 \bar{\eta}^{n+1} , \partial_1^2\bar{\eta}^{n} ,\nabla  \bar{q}^{n} )\|_{L^\infty_TH^2} . $$

Thus we immediately see from the above five estimates that
$$\{ ( u^n, \eta^n,
q^n, \eta^n_t,u^{n}_t)\}_{n=1}^\infty$$ is a Cauchy sequence in
$ C^0(\overline{I_T}, H^{ 3}_{\mm{s}}) \times C^0(\overline{I_T}, H^{ 3}_{\mm{s}})\times
  C^0(\overline{I_T},\underline{H}^3)\times C^0(\overline{I_T}, H^{2}_{\mm{s}} )\times C^0(\overline{I_T}, H^{2}_{\mm{s}} )$.
Hence there exists one limit function $(u^\kappa,\eta^\kappa,q^\kappa)$, which is the unique local solution of  the linear  $\kappa$-approximate problem \eqref{01dsaf16asdfasf0safd000} with $(u^\kappa,\eta^\kappa,q^\kappa)$ in place of $(\eta ,u ,q)$, and  also the unique local solution of  the linear  $\kappa$-approximate problem \eqref{01dsaf16asdfasf0000xx}.
In addition, it is easy to see that the solution $(u^\kappa,\eta^\kappa,q^\kappa)\in \mathbb{S}_{T  }$ by further using  Propositions \ref{pro:parabolic} and \ref{pro:elliptic}.

Noting that  the local time $T  _1$ is independent of the initial data and  the local solution constructed above satisfies $(u,\eta)|_{t=T  }\in {H}^{ 3}_{\mm{s}}\times {H}^{1,3}_{\mm{s}} $, thus, if $ T <\alpha $, we can further extend the local solution to be a global solution defined on $\overline{I_\alpha }$ by finite steps; moreover the obtained global solution is the unique solution of of  the linear  $\kappa$-approximate problem \eqref{01dsaf16asdfasf0000xx} and  belongs to $\mathbb{S}_{ \alpha  }$.
This completes the proof  of Proposition \ref{202102151544}.
\hfill$\Box$
\end{pf}

\subsection{Solvability of linearized  problem \eqref{01dsaf16asdfasf0000xx}}\label{subsubsce:03}

To investigate the solvability of linearized  problem \eqref{01dsaf16asdfsafasf0000}, we shall first $\kappa$-independent estimates of the solutions of the $\kappa$-approximate problem.

\begin{lem}\label{pro:0812}
Under the assumptions of Proposition \ref{202102151544},
we further assume that $\kappa \in (0,1)$, $\|\eta^0\|_{\underline{1},3}\leqslant \delta \in \mathbb{R}^+$, $(\partial_1 u^0, \partial_1^2 \eta^0) \in H^3$ and $\varsigma\in \mathbb{A}_{\alpha,\varsigma}^{3,1/4} $ satisfies
\begin{align}
& \partial_1\varsigma \in L^\infty_\alpha H^{1,3} ,\  \varsigma_t \in L^\infty_\alpha H^{1,3}\cap C^0(\overline{I_T}, H^2),\
\varsigma|_{t=0}=\eta^0 .   \nonumber
\end{align}
Then there exist polynomials
$P  ( I^0_H,\|({\varsigma},\partial_1{\varsigma},  {\varsigma}_t)\|_{L^{\infty}_{\alpha}H^{3}} )$,
$P  ( \|({\varsigma},\partial_1{\varsigma},  {\varsigma}_t)\|_{L^{\infty}_{\alpha}H^{1,3}} )$, positive constants $c$, $c_3\geqslant 1$,
and  $\delta_0>0$ (may depending on  $a$,  $m $ and $\Omega$),
such that, for any $\delta\leqslant \delta_0 $,
the local solution  $(\eta^\kappa,u^\kappa,q^\kappa)$ provided by Proposition \ref{202102151544} belongs to $\mathfrak{U}_{T_1}^{1,3}  \times \mathfrak{C}^0( \overline{I_{T_1}},{H}^{2,3}_{\mathrm{s}}) \times  C^0(\overline{I_{T_1}},\underline{H}^3)$
 and
enjoys the following estimates:
\begin{align}
&\label{2020081safdasaf21000}
\sup\nolimits_{t\in \overline{I_{T_1}}}\mathcal{E}_{\kappa}^L(t) +\kappa \| \partial_1^2\eta \|^2_{L^2_{T}H^{3}} \leqslant 4c_3^2{I^0_L}
\\
&\label{2020081safdasaf210000415}
{ \sup\nolimits_{t\in \overline{I_{T_1}}}\mathcal{E}^H_{\kappa}(t) \leqslant 4c_3^2{I^0_H}  }  ,\\
&\label{2020081609300416}
   \| q^\kappa\|_{C^0(\overline{I_{T_1}}, {H}^{3} ) } +\|(u_t^\kappa, \eta_t^\kappa,\nabla_{\mathcal{B}}q^\kappa )\|_{L^\infty_{T_1}H^3}\lesssim
P\left( I^0_H ,\|({\varsigma}, {\varsigma}_t)\|_{L^{\infty}_{T_1}H^{3}}^2\right),
\end{align}
where $\mathcal{E}_{\kappa}^L(t):=\|(\eta^\kappa,\partial_1\eta^\kappa,u^\kappa)(t)\|_{3}^2$, $\mathcal{E}^H_{\kappa}(t):=\|(u^\kappa,\partial_1\eta^\kappa)(t)\|_{\underline{1},3}^2$ and
\begin{align}
T_1:=\min\{ 1/3c_5 P  ( \|({\varsigma},\partial_1{\varsigma},  {\varsigma}_t)\|_{L^{\infty}_{\alpha}H^{1,3}} ) ,\alpha, 1\}   .
\label{20220103021516}
\end{align}
Moreover, for $(\eta^\kappa,u^\kappa)$   restricted on $[0,T_1]$,
\begin{align}
&\partial_1 \nabla^3(u^\kappa,  \partial_1\eta^\kappa )
\mbox{ is uniformly continuous in }H^{-1}  . \label{202103021540}
\end{align}
\end{lem}
\begin{pf}Let $T\leqslant \min\{\alpha,1\}$,
\begin{align}
\|\eta^0\|_{\underline{1},3}\leqslant \delta\in (0,1].
\label{20220010504202237}
 \end{align}
 We define that
 $$\|\cdot \|_{*,i}^2:=\|(\cdot ,D^\tau_1\cdot) \|_{ i}^2 \mbox{ and }  \mathcal{E}^\tau_{\kappa}(t):=\|(u^\kappa,\partial_1\eta^\kappa)(t)\|_{*, 3}^2.$$ In addition, we denote  $(u^{\kappa},\eta^{\kappa},{q}^{\kappa})$
  by $(u, \eta,{q})$ for simplicity.
Next we establish the desired uniform estimates for $(u,\eta,q)$  by seven steps.

(1) \emph{Estimates of $\eta$}

Recalling \eqref{20200805100sadfsa5}, it is easy to see that
\begin{align}
&\sup\nolimits_{t\in\overline{I_T}}\| \eta (t)\|_{3}^2+\kappa^2\| \partial_1^2\eta \|^2_{L^2_{T}H^{3}} \lesssim
I^0_L +T\sup\nolimits_{t\in\overline{I_T}}\|u(t)\|_3^2.\label{202104101103}
\end{align}

By \eqref{20200805DSAFS100sadfsa5}, we have
\begin{align}
&\limsup_{s\to 0^+}\sup_{t\in\overline{I_s}}\| \eta (t)\|_{1,3}^2  \lesssim
 \| \eta^0\|_{1,3}^2.\label{20210410sdaf1103}
\end{align}
Since
\begin{equation}\label{parabsadfolic}
                              \begin{cases}
D^\tau_1\eta_t-\kappa D^\tau_1\partial_1^2\eta=D^\tau_1 w ,\\
D^\tau_1\eta|_{t=0}=\eta^0  ,\\
D^\tau_1 \eta_3|_{\partial\Omega} =0,
\end{cases}
\end{equation}
applying Proposition \ref{pro:parabolic} to \eqref{parabsadfolic} yields that
\begin{align}
&\limsup_{s\to 0^+}\sup_{t\in\overline{I_s}}\| \partial_1  D^\tau_1 \eta (t)\|_{3}^2  \lesssim
 \| \eta^0\|_{2,3}^2.\label{2021041sdaf0sdaf1103}
\end{align}

(2) \emph{Estimates of $q$}

Noting that $q$ satisfies
\begin{equation}\label{202104152037}
                              \begin{cases}
-\Delta_{{\mathcal{B}}}q=K^1&\mbox{in } \Omega,\\
{\mathcal{B}}_{3j}\partial_{i}q=0  &\mbox{on } \partial\Omega,
\end{cases}
\end{equation}
where $K^1$ is defined in \eqref{202102141605}. Applying the estimate \eqref{202008051safa442} to the above boundary-value problem, and then using \eqref{202008121505} and \eqref{202008121535}, we have
\begin{align}
\label{202008121055}
\|{q}\|_{C^0(\overline{I_T}H^3)}
&\lesssim P\left( \|{\varsigma}\|_{L^{\infty}_{T}H^3}\right)\|K^1\|_{L^{\infty}_{T}H^1} \lesssim P\left( \|({\varsigma},  {\varsigma}_t)\|_{L^{\infty}_{T}H^3}\right)\sup\nolimits_{t\in\overline{I_T}}
\|(u,\partial_1^2\eta)(t)\|_2.
\end{align}

(3) \emph{$L^2$-norm estimates of $(u,\partial_1\eta)$}.

Applying $\mathcal{X}$ to  \eqref{01dsaf16asdfasf0000xx}$_2$ yields that
\begin{align}\label{202104091638}
\partial_t\mathcal{X} u+a\mathcal{X} u=m \mathcal{X}\partial_1^{ 2}\eta-\mathcal{X}\nabla_{\mathcal{B}}q
\end{align}
Multiplying \eqref{202104091638} by $\mathcal{X} u$ in $L^2$   yields that (for a.e. $t\in I_T$)
\begin{align}\label{202008121145}
\frac{1}{2}\frac{\mm{d}}{\mm{d}t}\int|\mathcal{X} u|^2\mm{d}y
+\int(a |\mathcal{X}u|^2\mm{d}y-m (\mathcal{X} \partial_1^{ 2}\eta)\cdot \mathcal{X} u)\mm{d}y
=-\int (\mathcal{X}\nabla_{{\mathcal{B}}}{q})\cdot \mathcal{X}u\mm{d}y.
\end{align}
By \eqref{202008121505}, \eqref{202008121535n} and \eqref{202008121055}, we can estimate that
\begin{align}\label{202008121207}
-\int \mathcal{X}\nabla_{{\mathcal{B}}}{q}\cdot \mathcal{X} u\mm{d}y
 \lesssim
\|\mathcal{B}\|_{\underline{1},2}\|\nabla {q}\|_{\underline{1},0}\| u\|_{1,0}
 \lesssim P\left(  \|({\varsigma}, \partial_1{\varsigma}, {\varsigma}_t)\|_{L^{\infty}_{T}H^3}\right)
 \sup\nolimits_{t\in\overline{I_T}}\|(u, \partial_1^2\eta) \|_2^2.
\end{align}

In addition, by \eqref{01dsaf16asdfasf0000xx}$_1$ and the integration by parts, we get
\begin{align}\label{202008121210}
\begin{aligned}
-m   \int (\mathcal{X}\partial_1^{ 2}\eta)\cdot  \mathcal{X} u\mm{d}y
&=\frac{m   }{2}\frac{\mm{d}}{\mm{d}t}\int|\mathcal{X}\partial_1\eta|^2\mm{d}y
+m   \kappa\int|\mathcal{X}\partial_1^{2}\eta|^2\mm{d}y.
\end{aligned}
\end{align}
Plugging \eqref{202008121207}--\eqref{202008121210} into \eqref{202008121145} and then integrating the resulting inequality over $(0,t)$, we arrive at
\begin{align}
&
\|\mathcal{X} \left(u, \partial_1\eta\right)(t)\|_{0}^2
+{  \kappa}\int_0^t\|\mathcal{X} \partial_1^2\eta(s)\|_{0}^2\mm{d}s\nonumber \\
&\lesssim \| \left(u^0, \partial_1\eta^0\right)\|_{\underline{1},0}^2 + TP\left( \|({\varsigma}, \partial_1{\varsigma}, {\varsigma}_t)\|_{L^{\infty}_{T}H^3}\right)
\sup\nolimits_{t\in\overline{I_T}}  \|(u, \partial_1^2\eta)(t)\|_2^2.\label{202008121220}
\end{align}

(4) \emph{Curl estimates of $(u,\partial_1\eta)$.}

Applying $\mm{curl}_{{\mathcal{B}}}$ to \eqref{01dsaf16asdfasf0000xx}$_2$  yields that
\begin{align}\nonumber
\partial_t\mm{curl}_{{\mathcal{B}}}u+a\mm{curl}_{{\mathcal{B}}}u
=m  \mm{curl}_{{\mathcal{B}}} \partial_1^2\eta
+\mm{curl}_{ {{\mathcal{B}}}_t}u.
\end{align}
Let the multiindex $\alpha$ satisfy $|\alpha|\leqslant2$. Applying $\partial^{\alpha}\mathcal{X}$
to the above identity yields
\begin{align}\label{202005021542nnbm}
\begin{aligned}
\partial_t\partial^{\alpha} \mathcal{X} \mm{curl}_{{\mathcal{B}}}u
+a\partial^{\alpha}\mathcal{X} \mm{curl}_{{\mathcal{B}}}u
=&{m  }\partial_1\partial^{\alpha}\mathcal{X}\mm{curl}_{{\mathcal{B}}} \partial_1\eta
+K^{\mathcal{X} ,\alpha}_3 ,
\end{aligned}
\end{align}
where $K_3^{\mathcal{X},\alpha}:=
\partial^{\alpha}\mathcal{X}\mm{curl}_{{{\mathcal{B}}}_{t}} u-m  \partial^{\alpha}\mathcal{X}
\mm{curl}_{\partial_1{{\mathcal{B}}}} \partial_1\eta
$.

Multiplying \eqref{202005021542nnbm} by
$\partial^{\alpha}\mathcal{X}\mm{curl}_{ {\mathcal{B}}}u$ in $L^2$, we obtain
\begin{align}
&\frac{1}{2}\frac{\mm{d}}{\mm{d}t}\int|\partial^{\alpha}\mathcal{X}\mm{curl}_{{\mathcal{B}}}u|^2\mm{d}y
+a\int|\partial^{\alpha}\mathcal{X}\mm{curl}_{{\mathcal{B}}}u|^2\mm{d}y
\nonumber\\
&-m\int {  }\left(\partial_1\partial^{\alpha}\mathcal{X}\mm{curl}_{{\mathcal{B}}}\partial_1\eta\right)
\cdot\partial^{\alpha}\mathcal{X}\mm{curl}_{{\mathcal{B}}}u\mm{d}y =\int K_3^{\mathcal{X},\alpha}\cdot \partial^{\alpha}\mathcal{X}\mm{curl}_{{\mathcal{B}}} u\mm{d}y
=:I_8 .\label{2020081214246}
\end{align}
Using the integration by parts, \eqref{01dsaf16asdfasf0000xx}$_1$ and a regularity method, we obtain
\begin{align}
&-m  \int { }\left(\partial_1\partial^{\alpha}\mathcal{X}\mm{curl}_{{\mathcal{B}}}\partial_1\eta\right)
\cdot\partial^{\alpha}\mathcal{X}\mm{curl}_{{\mathcal{B}}}u\mm{d}y \nonumber \\
&=\frac{m }{2}\frac{\mm{d}}{\mm{d}t}
\int|\partial^{\alpha}\mathcal{X}\mm{curl}_{{\mathcal{B}}} \partial_1\eta |^2\mm{d}y
+{m  \kappa} \int|\partial^{\alpha}\mathcal{X}\mm{curl} \partial_1^2\eta |^2\mm{d}y
+ I_{9}+I_{10}, \label{202008121536}
\end{align}
where we have defined that
$$
\begin{aligned}
&I_{9}:=m \int(\partial^{\alpha} \mathcal{X}\mm{curl}_{{\mathcal{B}}} \partial_1\eta )\cdot
\partial^{\alpha}\mathcal{X}\left(\mm{curl}_{\partial_1{\mathcal{B}}}u
-\mm{curl}_{\partial_t{\mathcal{B}}} \partial_1\eta \right)\mm{d}y,\\
&I_{10}:=
 m  \kappa
\int\bigg(\left(\partial^{\alpha} \mathcal{X}\mm{curl}_{\partial_1{\mathcal{B}}}\partial_1\eta\right)
 \cdot\partial^{\alpha} \mathcal{X}\mm{curl}_{{\mathcal{B}}} \partial_1^2\eta
+ |\partial^{\alpha}\mathcal{X}\mm{curl}_{{\mathcal{B}}-I} \partial_1^2\eta |^2
 \\
&\qquad\quad + 2 \left(\partial^{\alpha}\mathcal{X}\mm{curl}_{{\mathcal{B}}-I}\partial_1^2\eta\right)
\cdot\partial^{\alpha}\mathcal{X}\mm{curl} \partial_1^2\eta \bigg)\mm{d}y.
\end{aligned}
$$

By the estimates \eqref{202008121505}, \eqref{202008121535} and \eqref{202008121535n}, we can estimate that,
\begin{itemize}
  \item for the case  $\mathcal{X}=\mm{Id}$,
\begin{align}
&\label{202104152110}
I_{8}+I_{9}
 \lesssim P(\|({\varsigma},\partial_1{\varsigma}, {\varsigma}_t)\|_{L^{\infty}_{T}H^3})
 \|(u,\partial_1\eta)\|_3^2,\\
&\label{202104152111}
I_{10}\lesssim \kappa \big( P(\| {\varsigma} \|_{L^{\infty}_{T}H^{1,3}})
\| \eta\|_{1,3}(\|\mm{curl}  \eta \|_{2,2}+  \|{\mathcal{B}}-I\|_2\| \eta\|_{2,3})
\nonumber \\
&\qquad \quad+
\|{\mathcal{B}}-I\|_2\left(\|{\mathcal{B}}-I\|_2\| \eta \|_{2,3}^2
+  \|  \eta \|_{2,3}\|\mm{curl}  \eta \|_{2,2}\right)\big).
\end{align}
  \item for the case  $\mathcal{X}=D^\tau_1$,
\begin{align}
&\label{202104152112}
I_{8}+I_{9}
 \lesssim P(\|({\varsigma},\partial_1{\varsigma}, {\varsigma}_t)\|_{L^{\infty}_{T}H^{1,3}})
 \|(u,\partial_1\eta)\|_{*,3}^2,\\
&\label{202104152114}
 I_{10} \lesssim  \kappa \big( P(\| {\varsigma} \|_{L^{\infty}_{T}H^{2,3}})
\|\partial_1\eta\|_{*,3} (\|\mm{curl}  \partial_1^2\eta \|_{*,2}+\|{\mathcal{B}}-I\|_{\underline{1},2}\| \partial_1^2\eta\|_{*,3})
\nonumber \\
&\qquad \quad+
\|{\mathcal{B}}-I\|_{\underline{1},2}(\|{\mathcal{B}}-I\|_{\underline{1},2}\|\partial_1^2 \eta \|_{*,3}^2
+  \|\partial_1^2 \eta \|_{*,3}\|\mm{curl} \partial_1^2\eta  \|_{*,2})\big).
\end{align}
\end{itemize}

Noting that $\varsigma|_{t=0}=\eta^0$, by Newton--Leibniz formula,   \eqref{202008121535nnn} with $t=0$, \eqref{202008121535}  and \eqref{20220010504202237}, we have
\begin{align}
\| ({\mathcal{B}}(t)-I)\|_{\underline{i},2}
&\lesssim
\left\| ({\mathcal{B}}^0-I)\right\|_{\underline{i},2}+
\int_0^t\left\| {\mathcal{B}}_{\tau}(s)\right\|_{\underline{i},2}\mm{d}s\nonumber\\
&\lesssim
\begin{cases}
 \|\eta^0\|_{3}
+T   P(\|({\varsigma}, {\varsigma}_t)\|_{L^{\infty}_{T}H^{3}})
& \mbox{for } i=0; \\
  \|\eta^0\|_{ {1},3}
+T   P(\|({\varsigma},   {\varsigma}_t)\|_{L^{\infty}_{T}H^{1,3}})
 & \mbox{for }i=1.
\end{cases}
\label{202103020944}
\end{align}
Making use of the above estimates \eqref{202104152110}--\eqref{202103020944}, Young's inequality and the fact ``$\kappa\in(0,1)$",   we deduce from \eqref{2020081214246}--\eqref{202008121536}  that
\begin{align}
& \frac{\mm{d}}{\mm{d}t}
 \|\mm{curl}_{ \mathcal{B}} ( u ,  \partial_1\eta)\|_{2}^2
+{c \kappa}\| \mm{curl} \eta \|_{2,2}^2\nonumber \\
&\lesssim P( \|({\varsigma},\partial_1{\varsigma}, {\varsigma}_t)\|_{L^{\infty}_{T}H^3})
 \|(u,\partial_1\eta)\|_3^2\nonumber
\\
&\quad + \kappa \left(\|\eta^0\|_{3}
+T {P}(\|({\varsigma}, {\varsigma}_t)\|_{L^{\infty}_{T}H^{3}})\right)
\| \eta\|_{2,3}^2 \label{202008121724}
\end{align}
and
\begin{align}\label{202104101500}
&\frac{1}{2}\frac{\mm{d}}{\mm{d}t}
 \|\mm{curl}_{{\mathcal{B}}} (u ,  \partial_1\eta )\|_{*,2}^2
+{c \kappa}\| \mm{curl} \partial_1^2\eta \|_{*,2}^2\nonumber \\
&\lesssim P( \|({\varsigma},\partial_1{\varsigma}, {\varsigma}_t)\|_{L^{\infty}_{T}H^{1,3}})
\|(u,\partial_1\eta)\|_{*,3}^2
\nonumber \\
&\quad +  \kappa  \left(\|\eta^0\|_{\underline{1},3}
+T  {P}(\|({\varsigma}, {\varsigma}_t)\|_{L^{\infty}_{T}H^{1,3}})\right)
\|\partial_1^2\eta\|_{*,3}^2.
\end{align}

In addition, similarly to \eqref{202103020944}, we have
\begin{align}
\|\mathcal{X}\mm{curl}f(t)\|_{ 2}&\lesssim\left \|\mathcal{X}
\left(\mm{curl}_{{\mathcal{B}}^0-I}f(t),
\mm{curl}_{{\mathcal{B}}}f(t),\mm{curl}_{\int_0^t {\mathcal{B}}_{s}\mm{d}s}f(t)\right)\right\|_{ 2}
\nonumber \\
&\lesssim \|\mathcal{X}\mm{curl}_{{\mathcal{B}}}f(t)\|_{ 2}\nonumber \\
&\quad +
\begin{cases}
 (\|\eta^0\|_{3}
+ T {P}(\|({\varsigma}, {\varsigma}_t)\|_{L^{\infty}_{T}H^{3}}))\|f(t)\|_{3}
 & \mbox{for } \mathcal{X}=\mm{Id}; \\
 (\|\eta^0\|_{\underline{1},3}
+T {P}(\|({\varsigma}, {\varsigma}_t)\|_{L^{\infty}_{T}H^{1,3}})\|f(t)\|_{*,3}
 &  \mbox{for }\mathcal{X}=D^\tau_1.
\end{cases}
\nonumber
\end{align}
Integrating \eqref{202008121724} and \eqref{202104101500} over $(0,t)$, resp.,
and then using \eqref{20220010504202237}, \eqref{20210410sdaf1103} and the above estimate, we conclude that
\begin{align}
&\|(\mm{curl}u,\mm{curl}\partial_1\eta)(t)\|_{ 2}^2
+{  \kappa}\int_0^t\| \mm{curl}  \eta (s)\|_{2,2}^2\mm{d}s\nonumber \\
&\lesssim    I^0_L +(\|\eta^0\|_3^2+T
P(  \|({\varsigma},\partial_1{\varsigma},  {\varsigma}_t)\|_{L^{\infty}_{T}H^3}))
\sup\nolimits_{t\in\overline{I_T}} \|(u,\partial_1\eta)(t)\|_3^2\nonumber \\
&\quad+ \kappa  \left(\|\eta^0\|_{3}
+T  {P}(\|({\varsigma}, {\varsigma}_t)\|_{L^{\infty}_{T}H^3})\right)
 \int_0^t\| \eta(s)\|_{2,3}^2\mm{d}s
\label{202008121745}
\end{align}
and
\begin{align}\label{202104101434}
&\|(\mm{curl}u,\mm{curl}\partial_1\eta)(t)\|_{*,2}^2
+{   \kappa}\int_0^t\| \mm{curl} \partial_1^2\eta (s)\|_{*,2}^2\mm{d}s\nonumber \\
&\lesssim  I^0_H
+( \|\eta^0\|_{\underline{1},3}+TP(\|({\varsigma},\partial_1{\varsigma},  {\varsigma}_t)\|_{L^{\infty}_{T}H^{1,3}}))
\sup\nolimits_{t\in\overline{I_T}}\ \|(u,\partial_1\eta)(t)\|_{*,3}^2\nonumber \\
&\quad+ \kappa P\left(\|\eta^0\|_{\underline{1},3}
+T   (\|({\varsigma}, {\varsigma}_t)\|_{L^{\infty}_{T}H^{1,3}})\right)
 \int_0^t\|\partial_1^2\eta(s)\|_{*,3}^2\mm{d}s.
\end{align}

(5) \emph{Divergence estimates of $(u,\partial_1\eta)$.}

Similarly to \eqref{202103020944},   we have
\begin{align}
\|\mathcal{X}\mm{div}f(t)\|_{ 2}\lesssim \left\|\mathcal{X}\left(\mm{div}_{{\mathcal{B}}^0-I}f(t),\mm{div}_{{\mathcal{B}}}f(t),\mm{div}_{\int_0^t  {\mathcal{B}}_{\tau} \mm{d}\tau } f(t)\right) \right\|_{ 2} . \nonumber
\end{align}
Noting that $\mm{div}_{{\mathcal{B}}}u=0$, thus taking $f=u$ in the above estimate yields
\begin{align}
\|\mathcal{X}\mm{div}u \|_{ 2}
\lesssim &
\begin{cases}
(\|\eta^0\|_{ 3}
+T {P}(\|({\varsigma},{\varsigma}_t)\|_{L^{\infty}_{T}H^3}))
\|u\|_{3}
& \mbox{for } \mathcal{X}=\mm{Id} ; \\
(\|\eta^0\|_{\underline{1},3}
+T {P}(\|({\varsigma}, {\varsigma}_t)\|_{L^{\infty}_{T}H^{1,3}}))
\|u\|_{*,3}
 & \mbox{for }\mathcal{X}=D^\tau_1.
\end{cases}\label{202008121802}
\end{align}

Applying $\mm{div}_{{\mathcal{B}}}\partial_1$ to \eqref{01dsaf16asdfasf0000xx}$_1$
and then use \eqref{01dsaf16asdfasf0000xx}$_3$ yields
\begin{align}
 \partial_t \mm{div}_{{\mathcal{B}}}\partial_1\eta
-\kappa\partial_1\mm{div}_{{\mathcal{B}}}\partial_1^2\eta =\mm{div}_{ {\mathcal{B}_t}}\partial_1\eta-\mm{div}_{\partial_1{\mathcal{B}}}u
-\kappa\mm{div}_{\partial_1{\mathcal{B}}}\partial_1^2\eta=:K^{ 4}. \nonumber
\end{align}
We further derive from the above identity that
\begin{align}
& \frac{1}{2} \frac{\mm{d}}{\mm{d}t}\| \partial^{\alpha }\mathcal{X}\mm{div}_{\mathcal{B}}
\partial_1\eta   \|^2_0 +
\kappa  \| \partial^\alpha\mathcal{X} \mm{div}_{{\mathcal{B}}}\partial_1^2\eta  \|_2^2  \nonumber \\
&= \int \big( \partial^{\alpha } \mathcal{X}K^4 \partial^\alpha\mathcal{X} \mm{div}_{{ \mathcal{B}}}\partial_1\eta
-\kappa \partial^\alpha\mathcal{X} \mm{div}_{{\mathcal{B}}}(\partial_1^2\eta \partial^\alpha\mathcal{X} \mm{div}_{{\partial_1 \mathcal{B}}}\partial_1\eta )\big)\mm{d}y .  \label{202104112141}
\end{align}
Following the arguments of \eqref{202008121745} and \eqref{202104101434} by further using \eqref{20210410sdaf1103} and \eqref{2021041sdaf0sdaf1103},  we can derive from
\eqref{202104112141} that
\begin{align}
&\|\mm{div} \eta(t)\|_{1,2}^2
+ {\kappa} \int_0^t\|\mm{div} \eta(s)\|_{2,2}^2\mm{d}s\nonumber \\
&\lesssim \| \eta^0\|_{1,3}
+(\|\eta^0\|_3+TP(\|({\varsigma},\partial_1{\varsigma}, {\varsigma}_t)\|_{L^{\infty}_{T}H^3})\sup\nolimits_{t\in\overline{I_T}} \|(u,\partial_1\eta)(t)\|_3^2
\nonumber \\
&\quad
+\kappa \left(\|\eta^0\|_3+T P(\|({\varsigma},\varsigma_t)\|_{L^{\infty}_{T}H^3})\right)
\int_0^t\| \eta(s)\|_{2,3}^2\mm{d}s\nonumber \\
&\quad
+\kappa   P(\| {\varsigma} \|_{L^{\infty}_{T}H^{1,3}})
\int_0^t \|  \eta(s)\|_{1,3} \|  \eta(s)\|_{2,3}\mm{d}s  \label{202008122100}
\end{align} and
\begin{align}
&\| \mm{div}\partial_1\eta(t)\|_{*,2}^2
+ {\kappa} \int_0^t\|\mm{div}\partial_1^2 \eta(s)\|_{*,2}^2\mm{d}s\nonumber \\
&\lesssim \| \partial_1\eta^0\|_{\underline{1},3} +(\|\eta^0\|_{\underline{1},3}+T
P(\|({\varsigma},\partial_1{\varsigma}, {\varsigma}_t)\|_{L^{\infty}_{T}H^{1,3}})\sup\nolimits_{t\in\overline{I_T}} \|(u,\partial_1\eta)(t)\|_{*,3}^2
\nonumber \\
&\quad
+\kappa  \left(\|\eta^0\|_{\underline{1},3}
+T P(\|({\varsigma}, \varsigma_t)\|_{L^{\infty}_{T}H^{1,3}})\right)
\int_0^t\|\partial_1^2\eta(s)\|_{*,3}^2\mm{d}s\nonumber \\
&\quad
+\kappa   P(\| {\varsigma} \|_{L^{\infty}_{T}H^{2,3}})
\int_0^t\|\partial_1 \eta(s)\|_{*,3}\|\partial_1^2\eta(s)\|_{*,3} \mm{d}s .\label{202104101601}
\end{align}

Consequently,  we immediately deduce from \eqref{202008121745}--\eqref{202008121802}, \eqref{202008122100} and \eqref{202104101601} that
\begin{align}
&\|(\mm{div}u,\mm{curl}u,\mm{div}\partial_1\eta, \mm{curl}\partial_1\eta)(t)\|_{2}^2
+ {\kappa}\int_0^t\|(\mm{div} \eta,
\mm{curl} \eta)(s)\|_{2,2}^2\mm{d}s\nonumber \\
&\lesssim  I^0_L + (\|\eta^0\|_3+TP(\|({\varsigma},\partial_1{\varsigma}, {\varsigma}_t)\|_{L^{\infty}_{T}H^3})\sup\nolimits_{t\in\overline{I_T}}
\|(u,\partial_1\eta)(t)\|_{3}^2\nonumber \\
&\quad+   \kappa \left(\|\eta^0\|_{3}+T P(\|({\varsigma}, {\varsigma}_t)\|_{L^{\infty}_{T}H^3})\right)
\int_0^t\|\eta(s)\|_{2,3}^2
\mm{d}s\nonumber   \\ &\quad+\kappa   P(\| {\varsigma} \|_{L^{\infty}_{T}H^{1,3}})
\int_0^t \|  \eta(s)\|_{1,3} \|  \eta(s)\|_{2,3}\mm{d}s \label{202008122105}
\end{align}
and
\begin{align}\label{202104101622}
&\|(\mm{div}u,\mm{curl}u,\mm{div}\partial_1\eta, \mm{curl}\partial_1\eta)(t)\|_{*,2}^2
+ {\kappa}\int_0^t\|(\mm{div}\partial_1^2\eta,
 \mm{curl}\partial_1^2\eta)(s)\|_{*,2}^2\mm{d}s \nonumber \\
&\lesssim  I^0_H +(\|\eta^0\|_{\underline{1},3}+TP(\|({\varsigma},\partial_1{\varsigma}, {\varsigma}_t)\|_{L^{\infty}_{T}H^{1,3}})\sup\nolimits_{t\in\overline{I_T}
}\|(u,\partial_1\eta)(t)\|_{*,3}^2\nonumber \\
&\quad+ \kappa \left(\|\eta^0\|_{\underline{1},3}+T P(\|({\varsigma} , {\varsigma}_t)\|_{L^{\infty}_{T}H^{1,3}})\right)
\int_0^t\| \partial_1^2\eta (s)\|_{*,3}^2
\mm{d}s\nonumber\\
&\quad
+\kappa   P(\| {\varsigma} \|_{L^{\infty}_{T}H^{2,3}})
\int_0^t\|\partial_1 \eta(s)\|_{*,3}\|\partial_1^2\eta(s)\|_{*,3} \mm{d}s .
\end{align}

(6) \emph{Summing up the estimates of $(u,\eta,{q})$.}

Thanks to the estimates \eqref{202104101103},
\eqref{202008121220},
\eqref{202008122105}, {\eqref{202104101622}},
 the Hodge-type elliptic estimate \eqref{08171537} and Young's inequality, we have, \emph{for sufficiently small $\delta$},
\begin{align}
&\sup\nolimits_{t\in\overline{I_T}}\mathcal{E}_{\kappa}^L(t)
+  \kappa
\|\partial_1^2\eta^\kappa \|_{L^2_{{T}}H^{3}}^2\nonumber \\
&\leqslant c_5I^0_L
+c_5TP( \|({\varsigma},\partial_1{\varsigma}, {\varsigma}_t)\|_{L^{\infty}_{T}H^{3}})
\bigg(\sup\nolimits_{t\in\overline{I_T}}\mathcal{E}_{\kappa}^L(t)
+  \kappa \|\partial_1^2\eta^\kappa \|_{L^2_{{T}}H^{3}}^2\bigg) \label{202008122137}
\end{align}
and
\begin{align}\label{202104101640}
&\sup\nolimits_{t\in\overline{I_T}}\mathcal{E}^{\tau}_{\kappa}(t)
+ \kappa
\int_0^t \|\partial_1^2\eta^\kappa (s )\|_{*,3}^2\mm{d}s\nonumber \\
&\leqslant c_5I^0_H +c_5TP ( \|({\varsigma},\partial_1{\varsigma}, {\varsigma}_t)\|_{L^{\infty}_{T}H^{1,3}})
\bigg(\sup\nolimits_{t\in\overline{I_T}}\mathcal{E}^{\tau}_{\kappa}(t)
+\ \kappa \int_0^t \| \partial_1^2\eta^\kappa (s )\|_{*,3}^2\mm{d}s\bigg).
\end{align}

Obviously, we can choose a new  polynomial, still denoted by $P  ( \|({\varsigma},\partial_1{\varsigma},  {\varsigma}_t)\|_{L^{\infty}_{\alpha}H^{1,3}} )$, which is not less than the two polynomials in \eqref{202008122137} and \eqref{202104101640}.
Now we use $c_5$ and the new polynomial $P$ to define $T_1$ by \eqref{20220103021516}. Taking $T={T}_1$ in \eqref{202008122137} and \eqref{202104101640},  thus we get \eqref{2020081safdasaf21000} and
\begin{align}
&\sup\nolimits_{t\in\overline{I_T}}\mathcal{E}^{\tau}_{\kappa}(t)+ \kappa
\int_0^t \|\partial_1^2\eta^\kappa (s )\|_{*,3}^2\mm{d}s \leqslant  4c_3^2{I^0_H}   , \nonumber
\end{align}
which, together with \eqref{202104141653asfda}, yields that
\begin{align}
&\mm{ess}\sup\nolimits_{t\in {I_T}}\mathcal{E}_{\kappa}^H(t)  +    \kappa
\|\partial_1^2\eta^\kappa \|_{L^2_{{T}}H^{1,3}}^2\leqslant  4c_3^2{I^0_H}   .
\label{202105072000}
\end{align}
We further derive from \eqref{01dsaf16asdfasf0safd000}$_1$, \eqref{01dsaf16asdfasf0safd000}$_2$,  \eqref{202105072000} and the fact $(u_3,\partial_1^2\eta_3, \mathcal{B}_{3j}\partial_j q)|_{\partial\Omega}=0$ that
\begin{align}
\eta\in  C^0(\overline{I_T}, H^{1,3})\mbox{ and }u_t\in  C^0(\overline{I_T}, H^{2}_{\mm{s}}).
\label{2022010500720046}
\end{align}

By \eqref{2020081safdasaf21000} and \eqref{202008121055}, we get
\begin{align}
&\label{202008160942}
\|q\|_{C^0(\overline{I_T},H^3)}^2\lesssim
P\left(I^0_H , \|({\varsigma},   {\varsigma}_t {)}\|_{L^{\infty}_{T}H^3}^2\right).
\end{align}
Applying   \eqref{202008051442} to \eqref{202104152037}, and then using \eqref{202105072000}, we have
\begin{align}
\label{2020081safas21055}
 \|\nabla_{{\mathcal{B}}}q\|_{L^\infty_TH^{3}}
&\lesssim P\left( \|{\varsigma}\|_{L^{\infty}_{T}H^3}\right)\|K^1\|_{L^{\infty}_{T}H^2} \lesssim P\left(I^0_H , \|({\varsigma},   {\varsigma}_t {)}\|_{L^{\infty}_{T}H^3}^2\right).
\end{align}

In addition, exploiting \eqref{202105072000} and \eqref{2020081safas21055}, we deduce from \eqref{01dsaf16asdfasf0000xx}$_1$ and \eqref{01dsaf16asdfasf0000xx}$_2$  that
\begin{align}\label{202008asf160946}
\|(u_t,\eta_t)\|_{L^\infty_TH^3}
\lesssim
P\left(I^0_H , \|({\varsigma}, {\varsigma}_t {)}\|_{L^{\infty}_{T}H^3}^2\right).
\end{align}
Putting the estimates \eqref{202008160942}--\eqref{202008asf160946} together yields \eqref{2020081609300416} immediately.

(7)  \emph{Estimates of uniform  continuity of $\partial_1 \nabla^3(u^\kappa,  \partial_1\eta^\kappa )$ in $H^{-1}$. }

Let the multiindex $\beta$ satisfy $ |\beta|=  3$.
Similarly to \eqref{202210041111912}, we easily deduce from \eqref{01dsaf16asdfasf0000xx}$_1$ and \eqref{01dsaf16asdfasf0000xx}$_2$ that, for any $\varphi\in H^1 $ and for any $t$, $r\in I_T$,
\begin{align}
&\int\partial^\beta D_1^\tau\partial_1 (\eta(y,t)- \eta(y,r))
\cdot\varphi\mm{d}y=-\int_r^t\int\partial^\beta D_1^\tau(  u + \kappa \partial_1^2 \eta) \cdot \partial_1\varphi\mm{d}y\mm{d}s,  \label{2022104191625}\\
&  \label{2022104191625xxx}
\int \partial^\beta D_1^\tau( u (y,t) - u(y,r))  \cdot \varphi\mm{dy}\nonumber \\
&= \int_r^t\int
 ( \partial^{\beta }\nabla_{\mathcal{B}}q \cdot D_1^{-\tau} \varphi
- m  \partial^\beta D_1^\tau\partial_1 \eta \cdot \partial_1 \varphi
-a\partial^\beta D_1^\tau u  \cdot  \varphi )
\mm{d}y\mm{d}s.
\end{align}

Exploiting the uniform estimates \eqref{202105072000} and  \eqref{2020081safas21055}, we easily deduce that
\begin{align}
& D_1^\tau\nabla^3(u ,  \partial_1\eta  )
\mbox{ is uniformly continuous in }H^{-1} . \label{20210302asfd1540}
\end{align}
In addition, it is easy to check that
\begin{align}
\|D_1^\tau\nabla^3(u ,  \partial_1\eta  )\|_{L^\infty_TL^2}\lesssim\| \partial_1\nabla^3(u ,  \partial_1\eta  )\|_{L^\infty_TL^2} . \label{2021030safda2asfd1540}
\end{align}
Making use of \eqref{20210302asfd1540}, \eqref{2021030safda2asfd1540}, \eqref{2022205050220333} and the regularity $(u,\partial_1\eta )\in  C^0(\overline{I_T}, H^{3})$, we arrive at
\begin{align}
& \partial_1\nabla^3(u ,  \partial_1\eta )\in C^0(\overline{I_T}, L^2_{ \mm{weak}}) , \label{202103safafsda02asfd1540}
\end{align}
which, together with \eqref{202105072000}, yields \eqref{2020081safdasaf210000415}.
Noting that \eqref{2022104191625} and \eqref{2022104191625xxx} also hold with $\partial_1$  in place of $D_1^\tau$, thus we further derive
 the assertion in \eqref{202103021540}.

Thanks to \eqref{2022010500720046}, \eqref{202008asf160946}, \eqref{202103safafsda02asfd1540} and  the fact $(u ,\eta ,q)  \in \mathbb{S}_{T_1 }$, we easily see that
\begin{align}
(u ,\eta ,q )\in \mathfrak{U}_{T_1}^{1,3}  \times \mathfrak{C}^0( \overline{I_{T_1}},{H}^{2,3}_{\mathrm{s}}) \times  C^0(\overline{I_{T_1}},\underline{H}^3) . \nonumber
\end{align}
This completes the proof.
\hfill $\Box$
\end{pf}

Thanks to Lemma \ref{pro:0812}, we easily get a unique local solution  to the linearized problem   \eqref{01dsaf16asdfsafasf0000}  by a compactness argument and a regularity argument.
\begin{pro}\label{thm08}
Let the assumptions of Lemma \ref{pro:0812} be satisfied, and $\delta_0$, $T_1$ be provided by Lemma \ref{pro:0812}. Then, for any $\delta\leqslant \delta_0$,
the linearized problem \eqref{01dsaf16asdfsafasf0000}
defined on $\Omega_{T_1}$  admits a unique solution
$(u^{{L}},\eta^{{L}}, {q}^{L} )\in\mathfrak{U}_{T_1}^{1,3}\times   \mathfrak{C}^0( {\overline{I_{T_1}}},{H}^{2,3}_{\mm{s}})\times C^0(\overline{I_{T_1}},\underline{H}^3)$; moreover the solution satisfies
\begin{align}
&\label{202104191930}
 \|\eta^{{L}}\|_{C^0(\overline{I_{T_1}},H^3)}
\leqslant\|\eta^0\|_3+2 c_3\sqrt{ I^0_L }T,\\
&\label{202008131221}
 \|( u^{{L}},\eta^{{L}}, \partial_1 \eta^{{L}})\|_{C^0(\overline{I_{T_1}},H^3)}\leqslant 2 c_3
\sqrt{I^0_L},\\
&\label{2020081312210416}
\sup\nolimits_{t\in\overline{I_{T_1}}}  \|(u^{{L}}, \partial_1 \eta^{{L}})\|_{\underline{1},3}\leqslant 2 c_3 \sqrt{I^0_H},\\
&\label{2020081609540416}
  \|{q}^{L}  \|_{C^0(\overline{I_{T_1}},H^3)}+\| (u_t^{{L}}  ,\nabla_{\mathcal{B}}q^{{L}}  )\|_{L^\infty_{T_1}H^3 }
\leqslant P\left(  I^0_H ,\|({\varsigma},  \varsigma_t)\|_{L^{\infty}_{T_1}H^{3}}\right).
\end{align}
\end{pro}
\begin{pf}
Let $\kappa\in (0,1)$ and $(u^\kappa,\eta^\kappa,q^\kappa)\in\mathbb{S}_{{T_1} }$ be the local solution stated as in Lemma \ref{pro:0812}. Thanks to the $\kappa$-independent estimates \eqref{2020081safdasaf21000}--\eqref{2020081609300416} and \eqref{202103021540}
satisfied  by $(u^\kappa,\eta^\kappa,q^\kappa)$, we can easily follow the compactness argument as in the proof of Proposition \ref{pro:parabolic} to obtain a limit function $(u^{{L}},\eta^{{L}},q^{{L}})$, which is the solution of the linearized problem \eqref{01dsaf16asdfsafasf0000} and satisfies\begin{itemize}
\item the estimates
\begin{align}
&\label{20200813122fdas1}
 \|( u^{{L}},\eta^{{L}}, \partial_1 \eta^{{L}})\|_{L^\infty_{T_1}H^3}\leqslant 2 c_3
\sqrt{I^0_L},\\
&\label{202008131221041sdfa6}
 \|(u^{{L}}, \partial_1 \eta^{{L}})\|_{L^\infty_{T_1}H^{1,3}}\leqslant 2 c_3 \sqrt{I^0_H},\\
&\label{2020081609540fdsa416}
   \| (u_t^{{L}}  ,q,\nabla_{\mathcal{B}}q^{{L}}  )\|_{L^\infty_{T_1}H^3 }
\leqslant P\left(  I^0_H ,\|({\varsigma},  \varsigma_t)\|_{L^{\infty}_{T_1}H^{3}}\right).
\end{align}  \item the additional regularity
 \begin{align}
&\label{20200813sdafa122fdas1}
u^{{L}}\in {C^0(I_{T_1}H^2_{\mm{s}})}  ,\ \nabla^3 \partial_1 u \in {C}^0( \overline{I_T} ,L^2_{ \mm{weak}}),\  \eta^{{L}}\in \mathfrak{C}^0( \overline{I_{T_1}},{H}^{2,3}_{\mm{s}} ).
\end{align}
                                             \end{itemize}
In addition,
$$(q)_{\Omega}=0\mbox{ and }u_t|_{\partial\Omega}=0\mbox{ for a.e. }t\in I_T.$$

Thanks to Proportion \ref{pro:transport} and the regularity of $(u,\eta,q)$, we further have
\begin{align}
u^{{L}}\in {C^0(I_{T_1}H^3_{\mm{s}})}.
\end{align}

By \eqref{01dsaf16asdfsafasf0000}$_2$ and \eqref{01dsaf16asdfsafasf0000}$_3$, $q$ satisfies
\begin{equation}\nonumber
                              \begin{cases}
 \Delta_{{\mathcal{B}}}q=K^4:=m\mm{div}_{{\mathcal{B}}}  \partial_1^2\eta +\mm{div}_{{\mathcal{B}}_t}  u&\mbox{in } \Omega,\\
{\mathcal{B}}_{3j}\partial_{i}q=0  &\mbox{on } \partial\Omega,
\end{cases}
\end{equation}
where $K^4\in {C^0(I_{T_1}H^1 )}$. Applying Proposition \ref{pro:elliptic} to the above boundary-value problem, we easily obtain
\begin{align}
\label{202008121055xx}
  q \in {C^0(\overline{I_{T_1}},\underline{H}^3)} .
\end{align}
Thus we further derive from \eqref{01dsaf16asdfsafasf0000}$_2$ that
\begin{align}
\label{2020081sdfa21055}
  u_t \in C^0(\overline{I_{T_1}},H^2_{\mm{s}}).
\end{align}
Collecting the regularity of $(u,\eta,q)$ in \eqref{202008131221041sdfa6}--\eqref{2020081sdfa21055}, we easily see  that
 $(u^{{L}},\eta^{{L}},q^{{L}})\in \mathfrak{U}_{T_1}^{1,3}\times \mathfrak{C}^0( \overline{I_{T_1}},{H}^{2,3}_{\mm{s}} )\times  C(\overline{I_{T_1}},\underline{H}^3)$. Thus we further get \eqref{202008131221}--\eqref{2020081609540416} from \eqref{20200813122fdas1}--\eqref{2020081609540fdsa416}.

The estimate \eqref{202104191930} is obvious by \eqref{01dsaf16asdfsafasf0000}$_1$
satisfied by $(\eta^{{L}},u^{{L}})$ and the estimate \eqref{202008131221}.
In addition, the uniqueness can be easily verified by a standard energy method.
The proof of Proposition \ref{thm08} is complete.
\hfill$\Box$
\end{pf}

\subsection{Proof of Proposition \ref{202102182115}}\label{subsce:03}
Now we are in the position  to the proof of Proposition  \ref{202102182115}.
Let   $(u^0,\eta^0)$ satisfy the assumptions in Proposition \ref{202102182115}, $\|\eta^0\|_{\underline{1},3}\leqslant \delta\leqslant \{\delta_0, \iota/2\}$ and
\begin{align}
\nonumber
{T_2}=  \min\{1/3c_3 P(4c_3(b+1) , \delta/2c_3(b+1),1 \},
\end{align}
where the positive constants $c_3\geqslant 1$  and $\delta_0$ are provided by Proposition \ref{thm08} with $\alpha={ T_2}$.
By Proposition \ref{thm08}
and Lemma \ref{pro:1221}, for any $T\leqslant T_2$, we can construct a solution sequence
$$\{( u^{n},\eta^{n}, Q^{n})\}_{n=1}^\infty\subset\mathfrak{U}_{T }^{1,3}\times   \mathfrak{C}^0 (\overline{I_{T }},{H}^{2,3}_{\mm{s}} )\times {C^0(\overline{I_{T_1}},\underline{H}^3)} $$  such that
\begin{enumerate}[(1)]
\item $(u^{1}, \eta^1,{q}^{1})=(\eta^0, u^0,0)$;
  \item $(u^{n+1},\eta^{n+1}, {q}^{n+1})$ satisfies
\begin{equation}\label{01dsaf16asdfasf00002}
                              \begin{cases}
 \eta^{n+1}_t=u^{n+1} ,\\
u^{n+1}_t+\nabla_{\mathcal{A}^{n}}{q}^{n+1}+ au^{n+1}=m  \partial_1^2\eta^{n+1} ,\\
\div_{{\mathcal{A}}^{n}}u^{n+1}=0  , \\
(u^{n+1},\eta^{n+1})|_{t=0}=(u^0,\eta^0)  , \\
(u^{n+1}_3,\eta^{n+1}_3)|_{\partial\Omega}=\mathcal{A}^{n}_{3j}\partial_jq^{n+1}|_{\partial\Omega} =0
\end{cases}
\end{equation}
for $n\geqslant 1$, where $\mathcal{A}^{n}=(\nabla(\eta^{n}+I))^{-\mm{T}}$.
\item $(u^{n},\eta^{n}, q^{n})$ enjoys the estimates, for $n\geqslant 1$,
\begin{align}
&
  \| \eta^{n}  \|_{C^0(\overline{I_{T_1}},H^3)}^2 \leqslant 2\delta\leqslant \iota\leqslant 1,
\label{20dsfa2008131747} \\
&\label{202008131747} \|(u^{n},\eta^{n},\partial_1\eta^{n}) \|_{C^0(\overline{I_{T_1}},H^3)}  \leqslant 2c_3\sqrt{I^0_L}\lesssim b+1 , \\
&\label{202104101412}{
\sup\nolimits_{t\in\overline{I_T}} \|(u^{n}, \partial_1\eta^{n})(t)\|_{\underline{1},3} \leqslant 2c_3\sqrt{I^0_H } \lesssim b+1,} \\
&\label{202008131740n0416}
 \|q^n \|_{C^0(\overline{I_T},H^3)} + \|(u_t^{n} ,\nabla_{\mathcal{A}^{n-1}}q^n)\|_{L^\infty_{T}H^3}
\lesssim P ( b ) .
\end{align}
\end{enumerate}

Similarly to \eqref{202008121505} and \eqref{202008121535},
by using \eqref{20dsfa2008131747}  and \eqref{202008131747}, we have
\begin{align}
&\label{202008141600}
\sup\nolimits_{t\in  {I_{T}}} \|(\mathcal{A}^{n}_t,\mathcal{A}^{n},\partial_1\mathcal{A}^{n},1/J^{n})\|_2
\lesssim P(b ),\\
& \sup\nolimits_{t\in  {I_{T}}} \| \mathcal{A}^{n}-I\|_2
\lesssim \|\eta^n\|_3.\label{202008141dsafa600}
\end{align}

Similarly to \eqref{2022104191625} and \eqref{2022104191625xxx}, by a regularity method, we can easily verify that, for any $\varphi\in H^1 $ and for any $t$, $r\in \overline{I_T}$, $(\eta^{n+1}, u^{n+1})$  satisfies \begin{align}
&\int\partial^\beta \partial_1^2 (\eta^{n+1}(t)- \eta^{n+1}(r))
\cdot\varphi\mm{d}y=-\int_r^t\int\partial^\beta\partial_1  u^{n+1}  \cdot \partial_1\varphi\mm{d}y\mm{d}s  ,  \nonumber\\
&  \int \partial^\beta \partial_1(u^{n+1} (t) -u^{n+1}(r))  \cdot \varphi\mm{dy}\nonumber \\
& = \int_r^t\int
 ( \partial^{ \beta}   \nabla_{\mathcal{A}^{n}}q^{n+1}
-m \partial^\beta\partial_1^2 \eta^{n+1}+a\partial^\beta  u^{n+1}   )\cdot \partial_1\varphi
\mm{d}y\mm{d}s. \nonumber
\end{align}
Thus we further derive from the above two identities  that, for $n\geqslant 1$,
\begin{align}
\partial_1\nabla^3\left(u^{n}, \partial_1 \eta^{n}
\right)
\mbox{ are uniformly continuous in }H^{-1}.
\label{20210302154safdsadf0}
\end{align}
Next we shall prove that
$\{(\eta^{n}, u^{n}, {q}^{n})\}_{m=1}^{\infty}$
is a Cauchy sequence.

From now on, we always assume $n\geqslant 2$. We define that
$$\bar{\eta}^{n+1}:={\eta}^{n+1}-{\eta}^{n},\quad
\bar{u}^{n+1}:={u}^{n+1}-{u}^{n},\quad
\bar{{q}}^{n+1}:={{q}}^{n+1}-{{q}}^{n}\mbox{ and }\bar{\mathcal{A}}^{n}:={\mathcal{A}}^{n}-{\mathcal{A}}^{n-1}.$$
Then it follows from \eqref{01dsaf16asdfasf00002} that
\begin{equation}\label{01dsaf16asdfasf00003}
                              \begin{cases}
\bar{\eta}^{n+1}_t=\bar{u}^{n+1} ,\\
\bar{u}^{n+1} _t+\nabla_{\mathcal{A}^{n}}\bar{{q}}^{n+1}+a \bar{u}^{n+1}=
m    \partial_1^2\bar{\eta}^{n+1}
-\nabla_{\bar{\mathcal{A}}^{n}}{{q}}^{n} ,\\
\div_{{\mathcal{A}}^{n}}\bar{u}^{n+1}=-\div_{{\bar{\mathcal{A}}}^{n}}{u}^{n}   , \\
(\bar{u}^{n+1},\bar{\eta}^{n+1} )|_{t=0}=(0,0) , \\
(\bar{u}^{n+1}_3,\bar{\eta}^{n+1}_3)|_{\partial\Omega}=0   .
\end{cases}
\end{equation}

By \eqref{01dsaf16asdfasf00003}$_1$ and \eqref{01dsaf16asdfasf00003}$_4$, we have
\begin{align}\label{202008140740}
\begin{aligned}
 \|\bar{\eta}^{n+1}(t)\|_{C^0(\overline{I_{T}},H^3)}  \leqslant \int_0^t\|\bar{u}^{n+1}(s)\|_3\mm{d}s
\leqslant  {T}\|\bar{u}^{n+1}\|_{L^{\infty}_{T}H^3}.
\end{aligned}
\end{align}
Similarly to \eqref{202008141600}, we can estimate that
\begin{align}\label{202008140824}
& \sup\nolimits_{t\in  {I_T}}\|\bar{\mathcal{A}}^{n} \|_2
\lesssim
 \|\bar{\eta}^{n}\|_{{L^{\infty}_{T}H^3}} ,\\
& \sup\nolimits_{t\in  {I_T}}\|\bar{\mathcal{A}}^{n} \|_{1,2}
\lesssim
 \| \bar{\eta}^{n}\|_{{L^{\infty}_{T}H^{1,3}}} , \label{20200814sd0824}\\
\label{2020081608}
&\sup\nolimits_{t\in  {I_T}}
\| \bar{\mathcal{A}}^{n}_t \|_2
\lesssim  P(b)
 \|(\bar{u}^{n},\bar{\eta}^{n})\|_{{L^{\infty}_{T}H^3}} .
\end{align}

Thanks to \eqref{01dsaf16asdfasf00003}$_2$, $\bar{{q}}^{n+1}$ satisfies
\begin{equation}\label{01dsaf16asdfasf00004}
  \begin{cases}
 \Delta_{\mathcal{A}^{n}}\bar{{q}}^{n+1}
= \mm{div}_{\mathcal{A}^{n}}\left( m \partial_1^2\bar{\eta}^{n+1}
-\nabla_{\bar{\mathcal{A}}^{n}}{{q}}^{n}
-\bar{u}^{n+1}_t-a\bar{u}^{n+1}\right)
& \hbox{in }\Omega, \\
\nabla_{\mathcal{A}^{n}}\bar{{q}}^{n+1}\cdot\vec{\mathbf{n}} =
- \nabla_{\bar{\mathcal{A}}^{n}} q^{n}\cdot\vec{\mathbf{n}}  & \hbox{on }\partial\Omega,
\end{cases}
\end{equation}
where $\Delta_{\mathcal{A}^{n}}:=\mm{div}_{\mathcal{A}^{n}}\nabla_{\mathcal{A}^{n}}$.
Let $\zeta^n =\eta^n+y$. $(\zeta^n)^{-1}$ denotes the inverse function of $\zeta^n$ with respect to the variable $y$. We define that \begin{align}
&K^5:=
-\nabla_{\bar{\mathcal{A}}^{n}}{q}^{n} ,\
K^6:=  m  \partial_1^2\bar{\eta}^{n+1}
- \bar{u}^{n+1}_t- a  \bar{u}^{n+1}  -\nabla_{\bar{\mathcal{A}}^{n}}{q}^{n} ,\nonumber \\
&(\beta,\tilde{K}^5,\tilde{K}^6):=(\bar{q}^{n+1}, {K}^5,K^6)|_{y=(\zeta^n)^{-1}(x,t)},\nonumber
 \end{align}then, by \eqref{01dsaf16asdfasf00004},
\begin{equation}\nonumber
  \begin{cases}
 \Delta\beta=\mm{div}\tilde{K}^6 &\mbox{in }\Omega,\\
 \nabla \beta \cdot\vec{\mathbf{n}}= \tilde{K}^5 \cdot\vec{\mathbf{n}}  &\hbox{on }\partial\Omega.
\end{cases}
\end{equation}

Applying the elliptic estimate \eqref{neumaasdfann1n} to the above boundary-value problem and then making use of  \eqref{2021sfa04031901} and \eqref{2022104101908}, we have, for a.e. $t\in I_T$,
\begin{align}
\|\bar{q}^{n+1} \|_3\lesssim & P(b) \|\beta \|_3
 \lesssim P(b)\left(\|\mm{div}\tilde{K}^6\|_1+ \| \tilde{K}^5\|_2\right)\nonumber \\
\lesssim & P(b)\left(\| \mm{div}_{\mathcal{A}^{n}}{K}^6\|_1+\|{K}^5\|_2\right).
\label{2sfa02008141000}
\end{align}
Noting that
$$\mm{div}_{\mathcal{A}^{n}}{K}^6=\mm{div}_{\mathcal{A}^{n}}\left(
m  \partial_1^2\bar{\eta}^{n+1}-\nabla_{\bar{\mathcal{A}}^{n}}{q}^{n}-a \bar{u}^{n+1}\right)
 +\div_{{\mathcal{A}}^{n}_t}\bar{u}^{n+1}+
\partial_t \div_{{\bar{\mathcal{A}}}^{n}}{u}^{n} ,$$
thus, making use of \eqref{202008131740n0416}, \eqref{202008141600}, \eqref{202008140824} and \eqref{2020081608}, we easily estimate that
$$ \| \mm{div}_{\mathcal{A}^{n}}{K}^6\|_1+\|{K}^5\|_2
\lesssim P(b) \|\left(\bar{u}^{n}, \bar{u}^{n+1}, \bar{\eta}^{n}, \partial_1\bar{\eta}^{n+1}\right)\|_{{L^{\infty}_{T}H^3}} .$$
Putting the above estimate into \eqref{2sfa02008141000} yields
\begin{align}\label{202008141000}
\| \bar{q}^{n+1}\|_{C^0(\overline{I_{T}},H^3)}
 \lesssim P(b) \|\left(\bar{u}^{n}, \bar{u}^{n+1}, \bar{\eta}^{n}, \partial_1\bar{\eta}^{n+1}\right)\|_{{L^{\infty}_{T}H^3}} .
\end{align}

Similarly to the derivation of \eqref{202008121220}, we  have that
\begin{align}
\|(\bar{u}^{n+1},\partial_1\bar{\eta}^{n+1})(t)\|_{C^0(\overline{I_{T}},L^2)}  \lesssim  P(b)T\|\left(\bar{u}^{n}, \bar{u}^{n+1}, \bar{\eta}^{n}, \partial_1\bar{\eta}^{n+1}\right)\|_{{L^{\infty}_{T}H^3}}.
 \label{202008141635}
\end{align}

By \eqref{01dsaf16asdfasf00003}$_1$--\eqref{01dsaf16asdfasf00003}$_2$, we have
\begin{align}
&
\partial_t \mm{div}_{{\mathcal{A}}^{n}}\partial_1\bar{\eta}^{n+1}
 =\mm{div}_{ {\mathcal{A}}^{n}_t}\partial_1\bar{\eta}^{n+1}
-\partial_1 \mm{div}_{\bar{{\mathcal{A}}}^{n}}{u}^{n} -\mm{div}_{\partial_1{{\mathcal{A}}}^{n}}\bar{u}^{n+1},\nonumber \\
&
\partial_t \mm{curl}_{{\mathcal{A}}^{n}}\bar{u}^{n+1}
+a\mm{curl}_{{\mathcal{A}}^{n}}   \bar{u}^{n+1}
 - {m  }
\partial_1 \mm{curl}_{\mathcal{A}^{n}} \partial_1\bar{\eta}^{n+1}
\nonumber \\
&=\mm{curl}_{ {\mathcal{A}}_t^{n}}\bar{u}^{n+1}
- {m }\mm{curl}_{\partial_1\mathcal{A}^{n}}\partial_1\bar{\eta}^{n+1}  +  \mm{curl}_{\bar{\mathcal{A}}^{n}}(\nabla_{\mathcal{A}^{n-1}}q^{n}). \nonumber
\end{align}
Following the argument  of \eqref{202008122105} by further making use of \eqref{202104101412}, \eqref{202008141600}, \eqref{202008141dsafa600}, \eqref{202008140824} and \eqref{20200814sd0824}, we derive can from \eqref{01dsaf16asdfasf00003}$_3$ and the above identities that
\begin{align}
& \|(\mm{div}\bar{u}^{n+1},\mm{curl}\bar{u}^{n+1},\mm{div}\partial_1\bar{\eta}^{n+1},
  \mm{curl}\partial_1\bar{\eta}^{n+1})\|_{C^0(\overline{I_{T}},H^3)} \nonumber \\
&\lesssim  P(b)T
\|\left(  \bar{u}^{n+1},\bar{\eta}^{n}, \partial_1\bar{\eta}^{n}, \partial_1\bar{\eta}^{n+1}\right)\|_{{L^{\infty}_{T}H^3}}+  \|\eta^0\|_3\|\left( \bar{u}^{n+1}, \partial_1\bar{\eta}^{n+1}\right)\|_{{L^{\infty}_{T}H^3}}, \label{202008142035}
\end{align}
where we have used \eqref{202008140824} in the last inequality of \eqref{202008142035}.

Summing up the estimates \eqref{202008140740}, \eqref{202008141635} and \eqref{202008142035} and then using Hodge-type elliptic estimate \eqref{08171537}, we obtain,\emph{ for sufficiently small $\delta$,}
\begin{align}\nonumber
 \|\left(\bar{u}^{n+1},\bar{\eta}^{n+1}, \partial_1\bar{\eta}^{n+1} \right)\|_{C^0(\overline{I_{T}},H^3)}
 \lesssim P(b)T \|\left( \bar{u}^{n},\bar{u}^{n+1},\bar{\eta}^{n}, \partial_1\bar{\eta}^{n}, \partial_1\bar{\eta}^{n+1}\right)\|_{{L^{\infty}_{T}H^3}}
.
\end{align}
In particular
\begin{align}\label{202008142102}
 \|\left(\bar{u}^{n+1},\bar{\eta}^{n+1}, \partial_1\bar{\eta}^{n+1} \right)\|_{C^0(\overline{I_{T}},H^3)}
 \leqslant \|\left( \bar{u}^{n}, \bar{\eta}^{n}, \partial_1\bar{\eta}^{n} \right)\|_{{L^{\infty}_{T}H^3}}/2
\end{align}  \emph{for sufficiently small $T\in (0,T_2]$,} where the smallness depends  on $b$, $a$, $m $ and $\Omega$.

In addition, by \eqref{01dsaf16asdfasf00003}$_2$  and \eqref{202008141600}, we get that
\begin{align}\label{202008151041}
\|  \bar{u}^{n+1}_t\|_{C^0(\overline{I_{T}},H^2)}  \lesssim P(b)
 \|\left(\bar{u}^{n+1}, \bar{\eta}^{n}, \partial_1\bar{\eta}^{n+1},\bar{{q}}^{n+1}\right)\|_{{L^{\infty}_{T}H^3}}.
\end{align}

We immediately see from \eqref{202008141000}, \eqref{202008142102} and \eqref{202008151041} that the sequence $\{( u^{n},\eta^{n},{q}^{n},u_t^n)\}_{n=1}^{\infty}$ is a Cauchy sequence in $C^0(\overline{I_{T}},  H^3_{\mm{s}}\times H^{1,3}_{\mm{s}}\times \underline{H}^3\times H^{2}_{\mm{s}} )$. Thus
$$(u^{n}, \eta^{n}, {q}^{n},u_t^n)\to (u,\eta,{q},u_t)
\mbox{ strongly in }C^0(\overline{I_{T}}, H^3_{\mm{s}}\times H^{1,3}_{\mm{s}}\times \underline{H}^3\times H^{2}_{\mm{s}}).$$
Thanks to the strong convergence of $(u^{n},\eta^{n},{q}^{n},u^n_t)$, we can pass to limits in \eqref{01dsaf16asdfasf00002}, and then get a  solution $(u,\eta,{q},u_t)$ to the problem \eqref{01dsaf16asdfasf}.

Finally, thanks to \eqref{202008131747}--\eqref{202008131740n0416} and \eqref{20210302154safdsadf0}, we easily follow the  compactness argument and the regularity argument in Proposition \ref{thm08} to further derive that
\begin{align}
&(u,\eta,{q})\in  \mathfrak{U}_{T }^{1,3}\times \mathfrak{C}^0 (\overline{I_{T }},{H}^{2,3} )\times  C^0(\overline{I_{T}},\underline{H}^3),\  \nabla_{\mathcal{A} }q \in {L^\infty_TH^3}, \nonumber \\
&\label{2020081508001n}
 \|({u}, {\eta} , \partial_1{\eta})\|_{C^0(\overline{I_{T_1}},H^3)} \leqslant 2 c_3 \sqrt{I^0_L } ,\\
&\label{2020081508001nn}
\sup\nolimits_{t\in \overline{I_T}}\|({u} , \partial_1{\eta} )\|_{\underline{1},3} \leqslant 2 c_3 \sqrt{I^0_H}.
\end{align}
Finally it is easy to check that  $(\eta,u,q)$ is the unique solution to the problem
\eqref{01dsaf16asdfasf}, provided that $\delta$ is sufficiently small.
This completes the proof of Proposition \ref{202102182115}.
\hfill $\Box$

\appendix
\section{Analytic tools}\label{sec:09}
\renewcommand\thesection{A}
This appendix is devoted to providing some mathematical results, which have been used in previous sections.
 It should be noted that $\Omega$ appearing in what follows is  still defined by \eqref{0101a} with $h>0$ and we will also use the simplified notations appearing in Section \ref{subsec:03}. In addition,   $a\lesssim b$ still denotes $a\leqslant cb$, however the positive constant $c$ depends on the parameters and domain in the lemma, in which $c$ appears.
\begin{lem}\label{201806171834}
\begin{enumerate}[(1)]
 \item  Embedding inequality (see \cite[4.12 Theorem]{ARAJJFF}): Let $D\subset \mathbb{R}^2$ be a domain satisfying the cone condition, then
\begin{align}
&\label{esmmdforinfty}\|f\|_{C^0(\overline{D})}= \|f\|_{L^\infty(D)}\lesssim\| f\|_{H^2(D)}.
\end{align}
\item
Product estimates (see Section 4.1 in \cite{JFJSNS}): Let $D\in \mathbb{R}^2$ be a domain satisfying the cone condition, and the functions $f$, $g$ are defined in $D$. Then
\begin{align}
\label{fgestims}
&
 \|fg\|_{H^i(D)}\lesssim    \begin{cases}
                      \|f\|_{H^1(D)}\|g\|_{H^1(D)} & \hbox{ for }i=0;  \\
  \|f\|_{H^i(D)}\|g\|_{H^2(D)} & \hbox{ for }0\leqslant i\leqslant 2.
                    \end{cases}
\end{align}
                    \end{enumerate}
\end{lem}

\begin{lem}\label{10220830}
Poincar\'e inequality (see \cite[Lemma 1.43]{NASII04}): Let $1\leqslant p<\infty$, and $D$ be a bounded Lipchitz domain in $\mathbb{R}^n$ for $n\geqslant 2$ or a finite interval in $\mathbb{R}$. Then for any $w\in W^{1,p}(D)$,
\begin{equation}
\nonumber
\|w\|_{L^p(D)}\lesssim \|\nabla w\|_{L^p(D)}^p+\left|\int_{D}w\mathrm{d}y\right|^p.
\end{equation}
\end{lem}
\begin{rem}\label{10220saf830p}
By Poincar\'e inequality, we have, for any given $i\geqslant  0$,
\begin{align}
&\label{202012241002}
\| w\|_{1,i}\lesssim \|w\|_{2,i}\mbox{ for any  } w \mbox{ satisfying }\partial_1w,\ \partial_1^2w\in H^i.
\end{align}
 \end{rem}

 \begin{lem}
\label{lem:08171527}
 Hodge-type elliptic estimate:
Let $r\geqslant1$, then it holds that
\begin{equation}\label{08171537}
\|u\|_{r}\lesssim\|u\|_0+\|(\mm{curl}u,\mm{div}u)\|_{r-1}+\|u_3|_{\partial\Omega }\|_{H^{r-1/2}(\partial\Omega )}
\end{equation}
 \end{lem}
\begin{pf}  Please refer to \cite[Proposition 5.1]{Coutand1} for the proof.
\hfill $\Box$
\end{pf}

{\begin{lem}\label{pro:1221}Diffeomorphism mapping theorem:
There exists a sufficiently small constant $\iota\in(0,1]$, depending on $\Omega$, such that,
for any ${\varsigma}\in H_{\mm{s}}^3$
satisfying  $\|{\varsigma}\|_3\leqslant \iota$,
$\psi:=\varsigma+y$ satisfies $\inf_{y\in \Omega}\det(\nabla \varsigma  +I)\geqslant 1/4$
and the diffeomorphism  properties
\begin{align}
& \psi|_{y_3=i}   : \mathbb{R}\to \mathbb{R}\mbox{ is a } C^1(\mathbb{R})\mbox{-diffeomorphism mapping for }i=0,\ 1,\nonumber \\
&\nonumber
\psi:={\varsigma}+y:\overline{\Omega}\to\overline{\Omega} \mbox{ is a }C^1\mbox{-diffeomorphism mapping} .
\end{align}
\end{lem}
 \begin{pf} Please refer to \cite[Lemma A.10]{ORTIBMFIMVD} for the proof.
\hfill $\Box$
\end{pf}

\begin{lem}\label{20021032019018} Some results for functions with values in Banach spaces:
  Let $T>0$,  integers $i$, $j\geqslant 1$ be given and  $1\leqslant p\leqslant \infty$.
\begin{enumerate}
\item[(1)] Assume $f\in L^p_TH^{i}$, $\partial_1^k\nabla^{i}f\in L^p_TL^2$ for any $1\leqslant k\leqslant j$, then $f\in L^p_TH^{j,i}$.
\item[(2)] Let $X$ be a separable Banach space and $T>0$.
If $w\in C^0(\overline{I_T},X)$, then  $w$: $I_T\to X$ is a strongly measurable function and
\begin{align}
\|w\|_X\in C^0(\overline{I_T}).  \label{2022104121338}
\end{align}
\item[(3)] Let $  j\geqslant i+1$ and $c$ be a constant. If $f\in L^p(I_T,H^i)$ with $1<p\leqslant \infty$,  $\|f\|_{j}\leqslant c \|g\|_j$ holds for a.e. $t\in I_T$ and $g\in L^p(I_T,H^j)$, then
\begin{align}
f\in L^p_TH^i.
\label{202104032134}
\end{align}
\item[(4)]   Assume $f\in L^p_TH^{i}$ and $\{D_1^{\tau}f\}_{|\tau|\in (0,1)}$ is uniformly bounded in $L^p_TH^{i}$.
\begin{enumerate}
     \item[(a)]  There exists a sequence (still denoted by $D^\tau_1f$) of $\{D^\tau_1f\}_{|\tau|\in (0,1)}$ such that, if  $\tau\to 0$, then  \begin{align}
D_1^\tau f \rightharpoonup \partial_1f
\begin{cases}
\mbox{weakly-* in }L^\infty_TH^i&\mbox{for }p=\infty,\\
\mbox{weakly  in }L^p_TH^i &\mbox{for }1< p< \infty.\label{202104141653asfda}
\end{cases}
\end{align}
Moreover, $f \in L^\infty_TH^{1,i}$.
     \item[(b)] Let $p=\infty$.  If additionally $ \partial^\alpha f\in C^0(\overline{I_T},L^2_{\mm{weak}})$ and  $D_1^\tau\partial^\alpha f$ is uniformly continuous in $H^{-1}$, where  the multiindex $\alpha$ satisfying $|\alpha|=i$, then
\begin{align}
\label{2022205050220333}
D_1^\tau\partial^\alpha f \to \partial^\alpha \partial_1 f\mbox{ in } C^0(\overline{I_T},L^2_{\mm{weak}}) \mbox{ (some sequence)}.
\end{align}
  \item[(c)]  We additionally assume $f\in L^p_TH^{i}$  and $\{D^{-\tau}_1D^\tau_1f\}_{|\tau|\in (0,1)}$  is uniformly bounded in $L^p_TH^{i}$ with $1<p<\infty$, then (some sequence)
   \begin{align}D_1^{-\tau}(D_1^\tau f ) \rightharpoonup \partial_1^2 f \mbox{ weakly in }L^p_TH^i \mbox{ (some sequence)}.
\label{202104141653} \end{align}
\end{enumerate}
\item[(5)] Let $\varphi =\varsigma+y$ and $\varsigma\in \mathbb{A}_T^{3,1/4}$  defined by \eqref{2022104161633}.
\begin{enumerate}
\item[(a)] If $ C^0(\overline{I_T},H^i)$ or $f\in L^2_TH^i$ with $0\leqslant i\leqslant 3$,  then
 \begin{align}
\label{2020103250855}
&F:=f(\varphi ,t)\in C^0(\overline{I_T},H^i)\mbox{ or }L^p_TH^i
\end{align}
and
\begin{align}
 \mathcal{F}:=f(\varphi^{-1},t)\in C^0(\overline{I_T},H^i)\mbox{ or } L^p_TH^i.
\label{202104031901} \end{align}
Moreover,
\begin{align}
 & \|F\|_{L^p_TH^i}\lesssim
P(\|\varsigma\|_{L^\infty_TH^3}) \|f \|_{L^p_TH^i},
\label{2021sfa04031901} \\
&  \|\mathcal{F}\|_{L^p_TH^i}\lesssim
P(\|\varsigma\|_{L^\infty_TH^3})\|f \|_{L^p_TH^i}. \label{2022104101908} \end{align}
 \item[(b)]  If $ \varsigma $ additionally satisfies $\varsigma_t\in L^\infty_TH^2$, then for any $f$ satisfying $f\in L^p(I_T,H^i)$ and $f_t\in L^p(I_T,H^{i-1})$, where $1\leqslant i\leqslant 3$,  then \begin{align}F_t=(f_t(x,t) + \varsigma_t\cdot \nabla f(y,t))|_{x=\varphi}\in  L^p_TH^{i-1}\nonumber
\end{align}
and
\begin{align}\mathcal{F}_t=(f_t(y,t) - (\nabla \varphi)^{-1}\varsigma_t \cdot \nabla f(y,t))|_{y=\varphi^{-1}}\in  L^p_TH^{i-1}. \nonumber
\end{align}
 \end{enumerate}
\end{enumerate}
\end{lem}
\begin{pf} Please refer to \cite[Lemma A.10]{ORTIBMFIMVD} for the proof.
\hfill $\Box$
\end{pf}
}
\begin{lem}\label{lem:08181945}Elliptic estimates for Neumann boundary-value problem:  Let  $i\geqslant 2$ be an integer.  We further assume $A$ is a $3\times 3$ matrix function, each element of which belongs to  ${W^{ i-1,\infty}}( {\Omega})$, and satisfies positivity condition, i.e., there exists a positive constant $\theta$ such that $(A \xi )\cdot \xi\geqslant \theta |\xi|$ for a.e. $y\in \Omega$ and all $\xi\in \mathbb{R}^3$. If $f_1\in H^{i-2}$ and $f_2\in H^{i-1} $ satisfy the following compatibility condition
\begin{align}\label{202008181950n}
\int f_1\mm{d}y+\int_{\partial\Omega }f_2 \vec{\mathbf{n}}_3\mm{d}y_{\mm{h}}=0,
\end{align}
where $\vec{\mf{n}}_3$ denotes the outward unit normal vector $\vec{\mathbf{n}}$   on the boundary $\partial\Omega$,
 there exists a unique strong solution $p\in \underline{H}^{ j } $ to
the Neumann boundary-value problem of elliptic equations:
\begin{align}\label{neumann}
\begin{cases}
 -\mm{div}\left(A\nabla p \right)=f_1 &\mbox{in } \Omega,\\
   A_{3k}\partial_k p  =f_2     &\mbox{on } \partial\Omega .
\end{cases}
\end{align}
 Moreover $p$ satisfies
\begin{align}\label{neumaasdfann1n}
\|p\|_{ i}\lesssim  (1+\|A\|_{W^{{ i-1},\infty} (\Omega)})(\| f_1\|_{ i-2}+\|  f_2\|_{{ i-1}}).
\end{align}
\end{lem}
{\begin{pf}  Please refer to \cite[Lemma A.7]{ORTIBMFIMVD} for the proof.
\hfill $\Box$
\end{pf}
}

\vspace{4mm} \noindent\textbf{Acknowledgements.}
The research of Fei Jiang was supported by NSFC (Grant Nos. 12022102) and the Natural Science Foundation of Fujian Province of China (2020J02013), and the research of Song Jiang by National Key R\&D Program (2020YFA0712200), National Key Project (GJXM92579), and
NSFC (Grant No. 11631008), the Sino-German Science Center (Grant No. GZ 1465) and the ISF-NSFC joint research program (Grant No. 11761141008).

\renewcommand\refname{References}
\renewenvironment{thebibliography}[1]{%
\section*{\refname}
\list{{\arabic{enumi}}}{\def\makelabel##1{\hss{##1}}\topsep=0mm
\parsep=0mm
\partopsep=0mm\itemsep=0mm
\labelsep=1ex\itemindent=0mm
\settowidth\labelwidth{\small[#1]}%
\leftmargin\labelwidth \advance\leftmargin\labelsep
\advance\leftmargin -\itemindent
\usecounter{enumi}}\small
\def\newblock{\ }
\sloppy\clubpenalty4000\widowpenalty4000
\sfcode`\.=1000\relax}{\endlist}
\bibliographystyle{model1b-num-names}

\end{document}